\documentclass[11pt,twoside]{amsart}
 
\usepackage{amsfonts}
\usepackage{epsfig,graphics}
\usepackage{amssymb}
\usepackage{amscd}

\newtheorem{theoreme}{Theorem}[section]

\newtheorem{proposition}[theoreme]{Proposition}
       
\newtheorem{definition}[theoreme]{Definition}

\newtheorem{lemme}[theoreme]{Lemma}

\newcommand{\rk}{\operatorname{rk}}

\newcommand{\ad}{\operatorname{ad}}

\newcommand{\GL}{\operatorname{GL}}

\newcommand{\End}{\operatorname{End}}

\newcommand{\Ad}{\operatorname{Ad}}
\newcommand{\Aut}{\operatorname{Aut}}
\newcommand{\Diff}{\operatorname{Diff}}

\newcommand{\Hol}{\mathcal{H}ol}

\newcommand{\BZ}{{\bf Z}}
\newcommand{\NN}{{\bf N}} 
\newcommand{\RR}{{\bf R}}

\newcommand{\hm}{{\widehat{M}}}
\newcommand{\hx}{{\hat{x}}}
\newcommand{\hy}{{\hat{y}}}
\newcommand{\hz}{{\hat{z}}}

\newcommand{\dds}{{\frac{d}{ds}}}

\newcommand{\liek}{{\mathfrak{k}}}

\newcommand{\lien}{{\mathfrak{n}}}

\newcommand{\lieg}{{\mathfrak{g}}}
\newcommand{\lieh}{{\mathfrak{h}}}
\newcommand{\liep}{{\mathfrak{p}}}
\newcommand{\lies}{{\mathfrak{s}}}

\newcommand{\liea}{{\mathfrak{a}}}
\newcommand{\liem}{{\mathfrak{m}}}

\newenvironment{preuve}{\medskip \noindent {\bf Proof: }}
   {$\diamondsuit$ }

\begin{document}

\title{Topology of automorphism groups of parabolic geometries}
\author{Charles Frances and Karin Melnick}
\thanks{This project was initiated during a visit of the second author as Professeur Invit\'ee in May 2016 to the Universit\'e de Strasbourg, whom the authors thank for this support. Melnick partially supported by NSF grant DMS 1255462.} 
\date{\today}
\maketitle

\begin{abstract} We prove for the automorphism group of an arbitrary parabolic geometry that the $C^0$ and $C^\infty$ topologies coincide, and the group admits the structure of a Lie group in this topology.  We further show that this automorphism group is closed in the homeomorphism group of the underlying manifold.
\end{abstract}

\section{Introduction}

It is well known that the automorphism group of a rigid geometric structure is a Lie group. 
In fact, as there are multiple notions of rigid geometric structures, such as $G$-structures of finite type, Gromov rigid geometric structures, or Cartan geometries, the property that the local automorphisms form a Lie pseudogroup is sometimes taken as an informal definition of rigidity for a geometric structure.

There remains, however, some ambiguity about the topology in which this transformation group is Lie.  It is a subgroup of $\mbox{Diff}(M)$, assuming the underlying structure is smooth, so one may ask whether it admits the structure of a Lie group in the $C^\infty$, $C^m$ for some positive integer $m$, or even the compact-open, topology.  A related interesting question is whether the automorphism group is closed in $\mbox{Homeo}(M)$.

Theorems of Ruh \cite{ruh.finite.type} and Sternberg \cite[Cor VII.4.2]{sternberg.diffgeom} state that, if $H$ is the automorphism group of a $G$-structure of finite type of order $m$, then $H$ is a Lie group in the $C^m$ topology on $\Diff^{m+1}(M)$.  Gromov proved a similar result in \cite[Cor 1.5.B]{gromov.rgs} for a smooth Gromov-$m$-rigid geometric structure.  
%For those differential-geometric structures which, in sufficient order of differentiability $k+1$, determine a Cartan geometry, 
%% For essentially all differential-geometric structures corresponding to a Cartan geometry, it is also known that $H$ is a Lie group in the $C^k$ topology for appropriate $k$ (see \cite[Thm I.3.2]{kobayashi.transf}).
In the case of a smooth Riemannian metric $(M,g)$, the results above yield a Lie group structure for the $C^1$-topology on the isometry group $\mbox{Isom}(M,g)$.  

The classical theorems of Myers and Steenrod \cite{myers.steenrod}, however, 
         say that in this Riemannian case the $C^0$ and $C^m$ topologies coincide on $\mbox{Isom}(M,g)$ for all $m$. 
 Nomizu \cite{nomizu.affine.transf} proved the same for the group of affine transformations of a connection (under an assumption of geodesic completeness, which can be removed).  The essence of the proof is that exponential
   coordinates locally convert affine transformations to linear maps, and a sequence of linear transformations converging $C^0$ automatically converges $C^{\infty}$.
   
%% There are similar early vintage results of a global nature, for holomorphic automorphisms of a bounded complex domain \cite{cartan.holo.bounded.domain} (FIX REF), or for automorphisms of an almost-complex structure on a compact manifold \cite{boothby.kobayashi.wang}.  

%% For several specific differential-geometric structures, the $C^0$, or the compact-open, topology on $H \cong \Aut(M)$ also makes it a Lie group.  

%% For isometries of $C^{k+1}$ Riemannian metrics, $k \geq 1$, for example, the classical theorems of Myers and Steenrod \cite{myers.steenrod} say that the $C^0$ and $C^k$ topologies coincide.  

This article is concerned with the topology of local automorphisms of parabolic geometries 
(see section \ref{sec.definitions.examples} below for the general definition).  
These form a rich class of differential-geometric structures which behave differently 
from Riemannian metrics in the sense that their automorphisms can have strong dynamics, so, for example, a convergent sequence of automorphisms need not limit to a homeomorphism.  Parabolic geometries do not determine a connection; without the exponential map, it is 
no longer clear that a $C^0$-limit of smooth automorphisms should be smooth.  

\subsection{Statement of main results}

We first briefly survey some results for specific parabolic geometries, which will be generalized by our main theorem.  We remark that the first two theorems below, of Ferrand and Schoen, are proved by geometric-analytic techniques that are quite specific to the structures in question.

\begin{itemize}
\item In the course of proving the Lichnerowicz Conjecture on Riemannian conformal automorphism groups, Ferrand showed, 
using techniques of quasiconformal analysis, 
that if a homeomorphism $f$ is a $C^0$ limit of smooth conformal maps, then $f$ is also 
smooth and conformal \cite{lf.conf.regularity, lf.noncompact}.  

\item  Schoen \cite{schoen.cr} reproved Ferrand's result above, and extended it to strictly pseudoconvex $CR$-structures.
%as well as strictly pseudoconvex $CR$-structure, and also that if a homeomorphism $f$ is a $C^0$ limit of smooth 
%automorphisms of such a structure, then $f$ is a smooth automorphism. 
His proof uses scalar curvature and the conformal Laplace operator in the conformal case, and the analogous Webster scalar curvature and pseudoconformal subelliptic operator in the CR setting. 

\item In \cite{frances.degenerescence}, the first author proved for conformal pseudo-Riemannian structures that if a sequence of smooth
local conformal transformations converges $C^0$, then it converges $C^{\infty}$.  His approach is very different from the analytic techniques of \cite{lf.conf.regularity} and \cite{schoen.cr}: he uses the Cartan connection associated to these structures and the dynamics of the action on null geodesics.
\end{itemize}

We prove a generalization of the results recounted above to local automorphisms of arbitrary parabolic geometries.  Parabolic geometries are a broad family of geometric structures which nonetheless admit an extensive general theory.  Well known examples include the conformal semi-Riemannian structures and strictly pseudoconvex CR structures mentioned above, as well as more general nondegenerate CR structures, projective structures, and so-called path geometries, which encode ODEs (see \cite{cap.slovak.book.vol1} for a comprehensive reference).
 Definitions \ref{defi.parabolic.geometry} and \ref{defi.automorphism} below explain precisely what is meant by 
 \emph{parabolic geometry} and \emph{automorphism}/\emph{automorphic immersion}.  An automorphic 
 immersion can be informally defined as a differentiable immersion $f : U \to M$, where $U \subset M$ is an open set, which preserves the Cartan  
 geometry ${\mathcal C}$ on $M$.  When $U=M$ and $f$ is also a diffeomorphism, then $f$ is said to be an automorphism of 
 $(M, {\mathcal C})$. The set of automorphisms is a group that will be denoted $\Aut(M,{\mathcal C})$. Our main results can then be stated as follows:

\begin{theoreme}
\label{thm.compactopen}
Let $(M, \mathcal{C})$ be a smooth parabolic geometry.  Let $f_k: U \to M$ be a sequence of 
automorphic immersions of $(M,{\mathcal C})$ 
converging in the $C^0$ topology on $U$ to a map $h$. Then $h$ is smooth and $f_k \rightarrow h$ also in the $C^{\infty}$ topology.
\end{theoreme}

%% Now consider $H = \Aut(M,\hm, \omega)$.  Note that any open set in the $C^0$ topology is open in the $C^\infty$, and that these topologies are both second-countable.  Let $U$ be a $C^\infty$-open neighborhood of $1_H$.  We claim that U contains a $C^0$-open neighborhood $1_H$. If not, then for a countable $C^0$ neighborhood basis $(V_i)$ with $\cap V_i = \{ 1_H \}$, there would be $f_i \in V_i \backslash U$.  Thus $f_i$ converges in the $C^0$ topology to $1_H$, but not $C^\infty$, contradicting theorem \ref{thm.compactopen} above.  It follows that $U$ is also $C^0$ open.  We obtain the corollary:  

In section \ref{sec.proofcorollary} we will also prove the following:

 \begin{theoreme}
\label{thm.lietransfogroup}
Let $(M,\mathcal{C})$ be a smooth parabolic geometry. Then $\Aut(M,{\mathcal C})$ is
 a Lie transformation group in the compact-open topology.  Moreover, $\Aut(M,{\mathcal C})$ is closed in 
 $\operatorname{Homeo}(M)$ for this topology.
\end{theoreme}

\subsection{Definitions}
\label{sec.definitions.examples}

  Parabolic geometries are most conveniently defined in terms of Cartan geometries.  Let $G$ be a Lie group with Lie algebra $\lieg$, and  $P < G$ a closed subgroup.  We will assume throughout the article that the pair $(G,P)$ is \emph{effective}, meaning $G$ acts faithfully on $G/P$.  A noneffective pair can always be replaced by an effective one, with the same quotient space $G/P$ (see \cite{sharpe}).
 % that all {\bf Cartan geometries are effective}.  By this, we mean that in the model space,   This requirement is generally harmless (MAKE PRECISE EQUIVALENCE PRINCIPLE).  

\begin{definition}
\label{defi.cartan.geometry}
  A \emph{Cartan geometry} ${\mathcal C}$ on a manifold $M$, with model space 
${X} =G/P$ comprises $(\hm,\omega)$, where $\pi: \hm \rightarrow M$ is a principal $P$-bundle, and $\omega$ is a $\lieg$-valued one-form on $\hm$ satisfying:
\begin{itemize}
\item For all $\hat{x} \in \hm$, $\omega_{\hat{x}} : T_{\hat{x}} \hm \rightarrow \lieg$ is a linear isomorphism.

\item For all $g \in P$, $R_g^* \omega = (\Ad g)^{-1} \circ \omega$, where $R_g$ denotes the right translation by $g$ on $\hm$.

\item For all $X \in \liep$, $\omega(X^\ddag) \equiv X$, where $X^\ddag(\hx) = \left. \dds \right|_0 \hx. e^{sX}$.
\end{itemize}
\end{definition}
The basic example of a Cartan geometry modeled on ${X} =G/P$ is the \emph{flat} geometry on $X$ comprising $(G, \omega_G)$, where 
$\omega_G$ is the Maurer-Cartan form. 
 %% This model is {\it flat}.  General  Cartan geometry modelled on  $X=G/P$ have a Cartan curvature, the vanishing
 %% of which is an obstruction to be locally isomorphic to $X$.
 
\begin{definition} 
\label{defi.parabolic.geometry}
A \emph{parabolic geometry} is a Cartan geometry modeled on ${X} =G/P$, where $G$ is a semisimple Lie group with finite center and without compact local factors, and $P < G$ is a parabolic subgroup.
%In this case,
%$(\lieg,P)$ will be called a \emph{parabolic pair}.
\end{definition}

Our notion of parabolic subgroup is the standard one, which will be recalled in Section \ref{sec.notation}.

Essentially all classical rigid geometric structures correspond to a canonical Cartan geometry.  The process of canonically associating a Cartan geometry is called the \emph{equivalence problem} for a given geometric structure (see \cite{sharpe} for examples).  Parabolic geometries admit a uniform solution of the equivalence problem, in which each corresponds to a type of ``filtered manifold'' (barring one exception, projective structures); see \cite[Sec 3.1]{cap.slovak.book.vol1}, \cite{tanaka.parabolic}.

\begin{definition}  
\label{defi.automorphism}
For $(M,\mathcal{C})$ a smooth Cartan geometry with $\mathcal{C} = (\hm, \omega)$, an \emph{automorphism} is $f \in \mbox{Diff}(M)$ 
which lifts to a bundle automorphism ${\hat f}$ of $\hm$ satisfying ${\hat f}^* \omega = \omega$. 
The group of automorphisms is denoted $\mbox{Aut}(M,{\mathcal C})$.
% < \mbox{Diff}(M)$.

For an open subset $U \subseteq M$, a smooth immersion $f : U \to M$ is an \emph{automorphic immersion} of $(M,{\mathcal C)}$ if it lifts to a bundle map
${\hat f}: \widehat{U} = \pi^{-1}(U) \to \hm$  satisfying ${\hat f}^* \omega = \left.\omega \right|_{\widehat{U}}$. 
\end{definition}

As $(G,P)$ is effective, the elements $f \in \mbox{Aut}(M, \mathcal{C})$ correspond bijectively to their lifts ${\hat f}$ to $\hm$ satisfying ${\hat f}^*\omega =\omega$, and similarly for automorphic immersions  (see \cite[Prop. 3.6]{me.frobenius}).

\subsection{Lie topology on the automorphism group}
\label{sec.lietopology}

For ${\mathcal C} = (\hm, \omega)$ a smooth Cartan geometry on $M$, the group $\Aut(M,{\mathcal C})$ can be endowed with the structure of a Lie transformation group as follows
(we refer to the definition in \cite[Chap. IV]{palais.lie} of \emph{Lie transformation group}). The Cartan connection defines a framing ${\mathcal P}$ of $\hm$, the pullback by $\omega$ of any
   basis in $\lieg$. The automorphisms of a framing form a Lie transformation group; more precisely:
   
   \begin{theoreme}(S. Kobayashi \cite[Thm I.3.2]{kobayashi.transf})
    \label{theo.kobayashi}
    Let $N$ be a smooth, connected manifold with a smooth framing ${\mathcal P}$.  
    \begin{enumerate}
    \item $\Aut({\mathcal P}) < \mbox{Diff}(N)$ admits the structure of a Lie transformation group.
     \item For $m=0, \ldots, \infty$, the $C^m$-topology on $\Aut({\mathcal P})$ coincides with the Lie topology. 
     \item A sequence $f_k \in \Aut({\mathcal P})$ converges in the Lie topology if and only if there exists $z \in N$ such that $f_k(z)$ converges in $N$.
     \end{enumerate}
   \end{theoreme}

Denote by $\widehat{\Aut}(M,{\mathcal C})$ the group of bundle automorphisms of $\hm$ preserving $\omega$.  This is a $C^\infty$-closed subgroup of $\Aut(\hm, {\mathcal P})$, so it is closed in the Lie topology and inherits the structure of a Lie transformation group. 
%this Lie topology again coincides with the $C^k$ topologies $k=0, \ldots, \infty$.  
The isomorphism $\widehat{\Aut}(M,{\mathcal C}) \cong \Aut(M, {\mathcal C})$ then provides the latter with the structure of a Lie group, in fact of a Lie transformation group of $M$. The underlying topology on $\Aut(M,{\mathcal C})$, the pullback of the $C^\infty$ topology on $\widehat{\Aut}(M,{\mathcal C})$, will henceforth be referred to as the {\it Lie topology}.  For $U \subset M$, the automorphic immersions defined on $U$ admit a similarly defined topology, which we will also call the Lie topology.
   
Recall that the Lie topology on $\Aut(M,{\mathcal C})$, as well as all $C^m$-topologies, are second countable.  
%   It is thus legitimate to characterize those topologies thanks to the convergence of sequences. 
A sequence $(f_k)$ of automorphic immersions of $(M, \mathcal{C})$ converges in the Lie topology if and only if the lifted sequence $({\hat f}_k)$ converges $C^\infty$. 
%smoothly to some map in $\widehat{\Aut}(M,{\mathcal C})$.
Thus if $(f_k)$ converges for the Lie topology to an automorphic immersion, then it does for the $C^{\infty}$-topology on $M$.  In cases where $\hm$ is a subbundle of the $r$-frames of $M$, and ${\hat f}_k$ are the corresponding natural lifts of $f_k$, then $C^\infty$ convergence of $(f_k)$ on $M$ conversely implies convergence in the Lie topology.  Such is the case for many parabolic geometries, but this property in general is unclear.
Our proofs will go via the Lie topology on $\Aut(M,{\mathcal C})$, thus showing that it coincides with all $C^m$-topologies, $m=0,\ldots,\infty$, and similarly for automorphic immersions of $(M,\mathcal{C})$.
 
\subsection{Structure of the Proof}

A sequence $(f_k)$ of automorphic immersions converging in the $C^0$ topology gives rise to a \emph{holonomy sequence} $(p_k)$ in $P$.  The action of $(p_k)$ on $G/P$ reflects many features of the action of $(f_k)$ on $M$.  Section \ref{sec.ideas} contains the definition of holonomy sequences and their equicontinuity properties relative to those of $(f_k)$.  In Section \ref{sec.trans.to.model}, we translate the problem to a statement about holonomy sequences on $G/P$.  The proof of this statement, Theorem \ref{thm.induction}, proceeds by induction on $\mbox{rk}_\RR G$.  The base case, $\mbox{rk}_\RR G = 1$, is recalled from \cite{frances.lfrank1} in Section \ref{sec.proof.rank1}.  The task for the remainder of the paper is, given a holonomy sequence $(p_k)$ not conforming to the conclusion of Theorem \ref{thm.induction}, to find an invariant lower-rank subvariety of $G/P$ on which $(p_k)$ exhibits the same behavior, thus contradicting the induction hypothesis.  Section \ref{sec.tools} develops tools for identifying such a lower-rank subvariety, corresponding to certain manipulations on the root spaces of $\lieg$.  In Section \ref{sec.induction}, we apply these tools to complete the induction step.

\section{Holonomy and equicontinuity with respect to segments}
\label{sec.ideas}
Let $(M,\mathcal{C})$ be a Cartan geometry modeled on ${X} =G/P$, not necessarily parabolic. 
%Let $(f_k)$ be a sequence of automorphic immersions. 
%as in Theorem \ref{thm.compactopen}.  

\begin{definition} %[Equicontinuity]
  \label{def.equicontinuity}
  %Let $x \in M$.
  %We assume that there exists $x \in M$ such that $(f_k(x))$ converges to  $y\in M$.  
  A sequence $f_k : U \rightarrow M$ of automorphic immersions of $(M,\mathcal{C})$ is \emph{equicontinuous at $x$} $\in U$
   if there exists $y \in M$ such that for any $x_k \rightarrow x$ in $U$, the sequence $f_k(x_k) \rightarrow y$.
 \end{definition}

 If $f_k : U \rightarrow M$ converges $C^0$, then $(f_k)$ is clearly equicontinuous at every point of $U$.  
 %for the $$
 The following theorem says that conversely, equicontinuity \emph{at a single point} implies local $C^0$-convergence, at least for parabolic geometries.

\begin{theoreme}
 \label{thm.punctualequicontinuity}
 Let $(M,\mathcal{C})$ be a smooth parabolic geometry and $(f_k)$ a sequence of automorphic immersions equicontinuous at $x \in M$.  Then there exists an open neighborhood $U$ of $x$ on which a subsequence 
  of $(f_k)$ converges $C^{\infty}$ to a smooth map $h$.
\end{theoreme}

Note that Theorem \ref{thm.punctualequicontinuity} implies Theorem \ref{thm.compactopen}.
 
 \subsection{Holonomy sequences}

Let $f_k: U \rightarrow M$ be a sequence of automorphic immersions of $(M,\mathcal{C})$ which is equicontinuous at $x \in U$, with lifts $\hat{f}_k: \widehat{U} \rightarrow \widehat{M}$.  Associated to $(f_k)$ is a holonomy sequence $(p_k)$ in $P$, whose behavior around the base point $o = [P] \in G/P$ reflects much of the local behavior of $f_k$ around $x$.

% Then there exists $y \in M$ such that for every sequence $x_k \rightarrow x$, the images $f_k(x_k) \rightarrow y$.
 
 \begin{definition}
\label{defi.holonomy}
  Let $x_k \rightarrow x$ in $U$. 
  %$x\in M$ and $x_k \rightarrow x$, such that the sequence $(f_k(x_k))$ is contained in a compact subset of $M$.
  A sequence $(p_k)$ of $P$ is a \emph{holonomy sequence of $(f_k)$ along $(x_k)$} when there exist  $\hx_k \in \pi^{-1}(x_k)$ such that 
  $\{ \hx_k \}_{k \in \NN}$ and $\{ \hat{y}_k \} = \{ \hat{f}_k(\hx_k).p_k^{-1} \}_{k \in \NN}$ are
  bounded in $\hm$.  A \emph{holonomy sequence of $(f_k)$ at $x$} is any holonomy sequence along some sequence $x_k \to x$.
  
 \end{definition}

 We will denote by $\Hol(x)$ the set of all holonomy sequences of $(f_k)$ at $x$. Equicontinuity of $(f_k)$ at $x$ ensures that $\mathcal{H}ol(x)$ is nonempty.  Indeed, given $y \in M$ such that $f_k(x) \rightarrow y$, choose any $\hat{x} \in \pi^{-1}(x)$ and $\hat{y} \in \pi^{-1}(y)$.  Then there exists a sequence $(p_k)$ in $P$
   such that $\hat{f}_k(\hat{x}).p_k^{-1} \rightarrow \hat{y}$, so $(p_k) \in \mathcal{H}ol(x)$.
 
 %Observe that if $(f_k)$ is equicontinuous at $x \in M$, it admits holonomy sequences along any sequence $(x_k)$ converging to $x$.

 \subsection{Equicontinuity with respect to segments}

%The rough strategy for proving  Theorem \ref{thm.punctualequicontinuity} is to show that the 
Equicontinuity  of a sequence $(f_k)$ 
at $x$ will have strong consequences on the local behavior of its  holonomy sequences around the base point  $o \in G/P$.  A useful notion to capture this local behavior is \emph{equicontinuity with respect to segments}. 
  An \emph{unparametrized segment} in $G/P$ is a set of the form $[\xi]=\{ e^{t \xi}.o \ | \ t \in [0,1] \}$, for some 
 $\xi \in \lieg$.  Remark that distinct $\xi, \eta \in \lieg$ may define the same unparametrized segment.
 
 We fix a Riemannian
 metric in a fixed neighborhood of $o$ in $X$, with respect to which we will measure the length of segments $[\xi]$ in this neighborhood, and denote the results by $L([\xi])$.
 %The length of $[\xi]$ with respect to our Riemannian metric will be denoted $L([\xi])$.
 %segments (denoted $L([\xi])$).
%segments and $[\xi]_x \subset M$.  

\begin{definition}
 A sequence $(p_k)$ in $P$ is \emph{equicontinuous with respect to segments} if when a  sequence of segments  $[\xi_k]$ satisfies 
 $L([\xi_k]) \to 0$, and $p_k.[\xi_k]=[\eta_k]$, 
%for a sequence $(\eta_k) $ of $\lieg$, 
then every cluster value of $(\eta_k)$ in $\lieg$ 
 is in $\liep$.
 %, then one has $\eta_{\infty} \in \liep$.
 %\end{enumerate}
 %onsidering a subsequence of $(p_k)$, 
 %there exist bounded sequences $(\xi_k)$ and $(\eta_k)$ 
 %in $\lieg$,  such that 
 %$L([\xi_k]) \to 0$ while $p_k.[\xi_k]=[\eta_k]$ and $L([\eta_k]) \geq c>0$.
 %\end{enumerate}
\end{definition}

Observe that the condition $L([\xi_k]) \to 0$, hence the very notion of equicontinuity with respect to segments,  does not depend 
on the choice of Riemannian metric, since any two are bi-Lipschitz equivalent in a neighborhood of $o$.

\subsection{Relation of equicontinuity and equicontinuity with respect to segments}
\label{sec.relation}
%% Equicontinuity with respect to segments is a useful notion for relating the equicontinuity 
%% property of a sequence of automorphic immersions $(f_k)$ with the properties of its holonomy sequences. 

\begin{proposition}
 \label{prop.equivalence}
 Let $(M,\mathcal{C})$ be a Cartan geometry and $f_k: U \rightarrow M$ a sequence of automorphic immersions of $(M,\mathcal{C})$.
%modeled on ${X} =G/P$, for which $(f_k)$ is a sequence of geometric immersions.  
If $(f_k)$ is equicontinuous at $x\in U$, then every holonomy sequence  $(p_k) \in \Hol(x)$ is equicontinuous with respect to segments.
\end{proposition}

 The proof will use the development of curves $\gamma: [0,1] \to \hm$, a notion which we now recall.   Given such a smooth curve 
 $\gamma$,
    the equation $\omega_G(\tilde{\gamma}^{\prime}(s))=\omega(\gamma'(s))$, where $\omega_G$ is the Maurer Cartan form of $G$, defines 
     an ODE on $G$.  The solution $\tilde{\gamma}$ such that $\tilde{\gamma}(0)=id$ will be called the \emph{development of $\gamma$}.
    
The Cartan connection also yields an \emph{exponential map} on $\widehat{M}$: 
 any $u$ in $\lieg$ defines the $\omega$-constant vector field  $U^\ddag$ on $\widehat{M}$ by $\omega(U^\ddag) \equiv u$; denote 
 $\{ \varphi^t_{U^\ddag} \}$ the corresponding local flow. Observe that whenever $u \in \mathfrak{p}$, the flow  
 $\{ \varphi^t_{U^\ddag} \}$ is globally defined and corresponds to right multiplication by $e^{tu}$ in
 the bundle $\hm$ (by the third axiom in Definition \ref{defi.cartan.geometry}).
%We call $\phi_u^t$ the local flow  generated on $\hat{M}$ by the field $\hat U$.  
%Given 
%let ${\mathcal W}_{\hx} \subset \lieg$ be  the set of vectors $u$ such that  $\phi_u^t$ is defined for   $t \in [0,1]$ at $\hx$.  Then one defines 
The exponential map at $\hx \in \widehat{M}$ is defined in a neighborhood ${\mathcal U}={\mathcal U}_{\hat{x}}$ of the origin in $\lieg$ by
  $$u \mapsto \exp({\hx},u):= \varphi_{U^\ddag}^1.\hx,$$
%  The inverse function theorem guarantees that
Shrinking ${\mathcal U}$ if necessary makes the exponential map at $\hx$ a diffeomorphism onto a neighborhood of $\hx$ in $\hm$.
   For $u \in \mathcal{U} \subset \lieg$, we will denote
%Let us point out that for  a vector $u \in {\mathcal U}$, there are two distinct notions of exponentiation.
the exponential of $u$ at $\hx$ in $M$ by $\exp(\hx,u)$, and the exponential in the Lie group $G$ by $e^u$.
%We hope that those two distinct notations should avoid any confusion. 

It is easy to see that whenever $\hat{f}: \hm \to \hm$ is the lift of an automorphic immersion of $M$, then 
 $$ \exp(\hx,u)=\exp(\hat{f}(\hx),u).$$
 The $P$-equivariance property of $\omega$ leads to a corresponding equivariance property for the exponential map for all $p \in P$
  \begin{equation}
  \label{equivariance.property}
  \exp({\hx},u).p^{-1}=\exp({\hx.p^{-1}},(\Ad p).u) 
  \end{equation}  
%This relation is easily established 
% for every $ u \in {\mathcal W}_{\hx}, \, \, p \in P.$

 Last, we recall the following crucial reparametrization  lemma.
 \begin{lemme}[\cite{fm.nilpconf}, Proposition 4.3]
  \label{lem.reparametrage}
  Let $\gamma, \alpha : [0,1] \rightarrow \hm$ be smooth curves, with $\gamma(0)=\alpha(0)$, and let $q :[0,1] \to P$ be a smooth map
   satisfying $q(0)=id$. 
  \begin{enumerate}
  \item{Assume that for the developments $\tilde{\gamma}$ and $\tilde{\alpha}$, the relation
    ${\tilde \gamma}(s)={\tilde \alpha}(s).q(s)$ holds in $G$ for every $s \in [0,1]$.  Then $\gamma(s)=\alpha(s).q(s)$ holds in $\hm$.}
  \item{In particular, if $u, v \in \lieg$, and if there exists a smooth $a: [0,1] \to [0,1]$, with $a(0)=0$ and $a(1)=1$,  such that  
  $$ e^{su}=e^{a(s)v} q(s) \qquad \forall s \in [0,1],$$
 then, for every $\hy \in \hm$  such that $\exp(\hy,u)$ or $\exp(\hy,v)$ is defined,
  $$ \exp(\hy,u)=\exp(\hy,v).q(1)$$}
  \end{enumerate}
 \end{lemme}

%The proof of Lemma \ref{lem.reparametrage} can be found in REFERENCES.
%to section \ref{sec.preuvelemme}, and proceed with the 

\begin{preuve}(of Proposition \ref{prop.equivalence})  Assume for a contradiction that $(f_k)$ is equicontinuous at $x$, but that some holonomy sequence 
$(p_k)$ of $(f_k)$ at $x$ does not act equicontinuously with respect
 to segments.   Then $\hat{y}_k = \hat{f}_k(\hx_k).p_k^{-1}$ is bounded for a bounded sequence $(\hx_k)$ projecting to $x_k \to x$.
  After passing to a subsequence, we can assume $\hx_k \to \hx$ and $\hat{y}_k \to \hy$.
  
 Since $(p_k)$ is not equicontinuous with respect to segments, passing again to a subsequence, there exists a sequence of segments $[\xi_k]$, with $L([\xi_k]) \to 0$, as well
  as a sequence $(\eta_k)$ in $\lieg$ converging to $\eta_{\infty} \not \in \liep$, such that for all $k$:
   \begin{equation}
   \label{eq.chemins}
   p_k.[\xi_k]=[\eta_k].
   \end{equation}
   This condition can be expressed by the relation, valid for all $s \in [0,1]$: 
    $$e^{s Ad(p_k)(\xi_k)}=e^{\varphi_k(s) \eta_k}.p_k(s).$$
    Here,  $p_k : [0,1] \to P$ denotes a smooth path with $p_k(0) = 0$ and $\varphi_k: [0,1] \to [0,1]$ a
     nondecreasing diffeomorphism.  Given $\lambda > 0$ arbitrary small, let $0 < \lambda_k < 1$ be such that $\varphi_k(\lambda_k)=\lambda$ for all $k$. Then write 
\begin{equation}
\label{eq.reparam}    
e^{s Ad(p_k)(\lambda_k \xi_k)}=
     e^{\frac{\varphi_k(\lambda_k s)}{\varphi_k(\lambda_k)} 
     \varphi_k(\lambda_k) \eta_k}.p_k(\lambda_k s).
\end{equation}

Note that $L([\lambda_k \xi_k]) \to 0$. 
%since $[\lambda_k \xi_k]$ is a subsegment of $[\xi_k]$.  It follows that 
Thus for $\lambda$ sufficiently small, we can 
       replace $\xi_k$ and $\eta_k$ by
  $\lambda_k \xi_k$ and $\varphi_k(\lambda_k) \eta_k$, so that (\ref{eq.chemins}) holds, with the extra property that 
  $\exp(\hat{y}_k,\eta_k)$ is defined for all $k \in \NN$, and $\eta_{\infty}$ is in an injectivity domain of the map
  $u \mapsto \exp(\hy,u)$. In 
  particular, if we call $y:=\pi(\hy)$, the fact that  $\eta_{\infty} \not \in \liep$ implies, shrinking $\lambda$ again if necessary, $\pi(\exp(\hy,\eta_{\infty})) \not =y$.

The next step is to show that $\pi(\exp(\hx_k,\xi_k))$ is defined for $k$ large enough, and converges to $x$.  To this aim, define a left-invariant Riemannian metric $\rho_G$ on $G$ by left translating
  any  scalar product  $< \ ,\ >$ on  $\lieg$, and a corresponding Riemannian metric $\rho$ on $\hm$, with
   $$ \rho(u,v):= \langle \omega(u),\omega(v) \rangle.$$
    By the definition of $\rho$, if $\gamma$ is a curve in $\hm$ and $\tilde{\gamma}$ its development in $G$, then
     $L_{\rho_G}(\tilde{\gamma})=L_{\rho}(\gamma)$.  Fix $\epsilon>0$ small enough that $\forall \ k \in \NN$, the $\rho$-ball 
     $B(\hx_k,\epsilon)$
      of center $\hx_k$ and radius $\epsilon$ has compact closure in $\hm$.

     Now consider  the curve $s \mapsto e^{s \xi_k}$. We fix  $\Sigma$ a small submanifold of $G$ containing $1_G$, which is transverse to the 
      fibers of $\pi_X: G \to X=G/P$, and such that the restriction of $\pi_X$ to $\Sigma$  yields a diffeomorphism  $\psi : \Sigma \to U$, where $U$ is
      a neighborhood  of $o$ in $X$.  
      For $k$ large enough, there exists a smooth $q_k : [0,1] \to P$, with $q_k(0)=id$, such that
      $\alpha_k(s)=e^{s \xi_k}.q_k(s)$ is contained in $\Sigma$.  Of course $\psi(\alpha_k([0,1]))=[\xi_k]$.
       %(***$[\alpha_k]$ doesn't fit the definition of unparametrized segment, so should clarify the notation here***).
        Two Riemannian metrics on $\Sigma$ are always locally bi-Lipschitz equivalent, hence there exist $C_1, C_2 > 0$
         such that for $k$ large enough:
         $$ C_1L([\xi_k]) \leq L_{\rho_G}(\alpha_k) \leq C_2 L([\xi_k]).$$
      We infer that $L_{\rho_G}(\alpha_k) \to 0$; in particular, for $k \geq k_0$, $L_{\rho_G}(\alpha_k) < \epsilon$.  Now consider,
      for each $k \geq k_0$, the 
      first-order ODE on $\widehat{M}$:
       
       \begin{equation}
        \label{eq.dev}
        \omega(\beta_k^{\prime})=\alpha_k^{\prime}
       \end{equation}

with initial condition $\beta_k(0)=\hx_k$.  
If $[0,\tau_k^*)$, is a maximal interval of definition for $s \mapsto \exp(\hx,s \xi_k)$, then for all $k$,
 $\beta_k(s):=\exp(\hx_k,s \xi_k).q_k(s)$, $s \in [0,\tau_k^*)$, is a maximal solution of our ODE, by Lemma \ref{lem.reparametrage}.  By the definition 
 of $L_{\rho}$, we have  
   $L_{\rho}(\beta_k)=L_{\rho_G}(\alpha_k)$.  If $\tau_k^* \leq 1$, the inequality  $L_{\rho_G}(\alpha_k) < \epsilon$ implies that $\beta_k$ is included in the relatively compact set 
  $B(\hx_k, \epsilon)$; this contradicts the maximality of $\tau_k^*$.  We thus infer $\tau_k^* >1$, which ensures that $\beta_k(1)$, hence 
  $\exp(\hx_k,\xi_k)=\beta_k(1).q_k(1)^{-1}$ is defined.  Moreover, $L_{\rho}(\beta_k)=L_{\rho_G}(\alpha_k) \rightarrow 0$, so $\beta_k(1) \rightarrow \hx$. Projecting to $M$ gives $\pi(\exp(\hx, \xi_k)) \to x$.

Now Lemma \ref{lem.reparametrage}, combined with equation (\ref{eq.reparam}) above says that for all $k \geq k_0$,
 $$ f_k(\exp(\hx_k,\xi_k).p_k^{-1})=\exp(\hat{y}_k,\Ad(p_k)\xi_k)=\exp(\hat{y}_k,\eta_k).p_k(1).$$  
 Projecting this relation on $M$, we obtain
 $$ \hat{f}_k(\pi(\exp(\hx_k,\xi_k)))=\pi(\exp(\hat{y}_k,\eta_k)).$$
 After possibly passing to a subsequence, the right-hand term converges to
 $\pi(\exp(\hy,\eta_{\infty})) \not = y$, while we just showed $\pi(\exp(\hx_k,\xi_k)) \to x$; this yields the desired contradiction with the equicontinuity of $(f_k)$ at $x$.
 \end{preuve}

%\subsection{Equicontinuity with respect to segments in the parabolic case}
% \subsection{Characterization of equicontinuity with respect to segments}

%  When the model pair $(\lieg,P)$ is parabolic, we obtain the following result:
%- Some notations about the roots, the group $A$. ($\Lambda$, the parabolic $P$? maybe later).

%Theorem \ref{thm.compactopen} will be  a straightforward consequence of 

%% \begin{theoreme}
%%  \label{thm.equicontinuity}
%%  Let $(M,\hm, \omega)$ be a Cartan geometry modeled on a parabolic pair $(\lieg,P)$.  Let $(f_k)$ be a sequence of automorphic immersions which is 
%%   equicontinuous at $x \in M$.  Then there exists a subsequence $(a_k) \in \Hol(x)$ such that 
%%   \begin{enumerate}
%%    \item $a_k \in A \ \forall k$.
%%    \item $\left. \Ad(a_k) \right|_{\lien^-}$ 
%%    is  bounded  in $\End(\lien^-)$.
%%   \end{enumerate}
%% \end{theoreme}

%Once point $1$ is proved, point $2$ is almost obvious, so the main difficulty is to prove the first point of Theorem \ref{thm.equicontinuity}.

%We will prove Theorem \ref{thm.equicontinuity} in the last section.  

\subsection{Vertical and transverse perturbations of holonomy sequences}
\label{sec.vertical.transverse}

Proposition \ref{prop.equivalence} translates equicontinuity of $(f_k)$ at $x$ to a property of sequences in $\Hol(x)$, which are in turn sequences of $P$ acting on $X = G/P$.  In this section we define several operations on sequences in $P$ which preserve $\Hol(x)$.

%The next step is to explore how equicontinuity with respect to segments restricts sequences of $\Hol(x)$.

%% The goal for the sequel is to prove Theorem \ref{thm.equicontinuity}.  Proposition \ref{prop.equivalence} provides the first step toward transferring the problem to the model space: if $p_k \in P$ form a holonomy sequence acting equicontinuously with respect to segments on $X = G/P$, then we wish to convert $(p_k)$ to a sequence of the desired form in such a way that the result is still a holonomy sequence for $(f_k)$ at $x$.

%  The set of reduced holonomy sequences is not completely random.  It is kind of structurated, and the aim of this paragraph is to 
%   exhibit some operations which leave it  invariant.

Holonomy sequences involve many choices: of $(x_k)$, of $(\hx_k)$, and of $(\hat{y}_k) = (\hat{f}(\hx_k) p_k^{-1})$, in the notation of Definition \ref{defi.holonomy}.  The right and left \emph{vertical perturbations} of $(p_k)$ correspond to other possible choices of $(\hx_k)$ and $(\hat{y}_k)$, respectively.  
 
\begin{definition}  
\label{defi.vertical.perturbation}
Let  $(p_k)$ be  a sequence in $P$. A {\it vertical perturbation of $(p_k)$} is a sequence $q_k=l_kp_km_k$ where 
  $(l_k)$ and $(m_k)$ are two bounded sequences in $P$.
   \end{definition}
  
%Observe that if $(f_k)$ is equicontinuous at $x \in M$, it admits holonomy sequences along any sequence $(x_k)$ converging to $x$.
\emph{Transverse perturbations} of $(p_k)$ correspond roughly to other possible choices of $(x_k)$ converging to $x$.

\begin{definition} 
\label{defi.transverse.perturbation}
For $(p_k)$  a sequence of $P$, a sequence $(q_k)$ of $P$ is said to be a \emph{transverse perturbation of $(p_k)$}  when there exist two
sequences $(\xi_k)$ and $(\eta_k)$ in $\lieg \backslash \liep$ such that:
\begin{enumerate}
\item{$ q_k=e^{- \eta_k}p_k e^{\xi_k}.$}
\item{The sequences $(\xi_k)$ and $(\eta_k)$ both converge to $0$.}
\item{For every $s \in \RR$, $e^{- s \eta_k}p_k e^{ s \xi_k}$ belongs to $P$.}
\end{enumerate}
\end{definition}

The other choice of $(x_k)$ in this case is $\pi(\exp(\hat{x}_k,\xi_k))$, as will be seen in the proof below.
 
 \begin{lemme}
   \label{lem.stability.transverse}
   Let $(M,\mathcal{C})$ be a Cartan geometry, and let $f_k : U \rightarrow M$ be a sequence
    of automorphic immersions.
   For any $x \in U$, the set of holonomy sequences $\Hol(x)$ is stable by vertical and transverse perturbations.  
  \end{lemme}

 %As we will see in the last section, the assumption $Nor(P^{\,o})=P$ is not really a serious one, and we can reduce our study to this case for any
  %Cartan parabolic geometry.
 
\begin{preuve}
 We consider $(p_k)$  a sequence belonging to $\Hol(x)$.  By definition, there exists $(\hx_k)$ a bounded sequence in $\hm$ 
 such 
  that $\hat{y}_k = \hat{f}_k(\hx_k).p_k^{-1}$ is bounded, and the projection $x_k$ on $M$ converges to $x$.
        
  Assume  that $(q_k)$ is obtained from $(p_k)$ by vertical perturbation, namely there exist  bounded sequences $(l_k)$ and $(m_k)$ in $P$ 
   such that $q_k=l_kp_km_k$.  Then $(\hx_k.m_k)$ is bounded in $\hm$, and still projects on $(x_k)$.
    Moreover 
$$\hat{f}_k(\hx_k.m_k)q_k^{-1}=\hat{y}_k.l_k^{-1}$$
 is still bounded in $\hm$.  It follows that $(q_k)$ is a holonomy sequence 
    at $x$.

  We now handle the  case of a transverse perturbation $q_k=e^{-\eta_k}p_ke^{\xi_k}$.   
%We have to show that $(q_k)$ is still a  holonomy sequence at $x$.
        The sequence $(\hx_k)$ is bounded and $\xi_k \rightarrow 0$, hence $(\hat{z}_k) = (\exp(\hx_k,\xi_k))$ is bounded in $\hm$, too; moreover, $\pi(\hat{z}_k)$ converges to $x$.  It remains to show that $\hat{f}_k(\hat{z}_k).q_k^{-1}$ is bounded  in $M$.
     Write this expression as $\hat{f}_k(\hat{z}_k).p_k^{-1}.p_kq_k^{-1} $.  
%The previous analysis allows to determine
By the equivariance (\ref{equivariance.property}) of the exponential map,
     $$\hat{f}_k(\hat{z}_k).p_k^{-1} = \exp(\hat{f}_k(\hx_k).p_k^{-1},\Ad(p_k)\xi_k).$$  
     
    Point $(2)$ in the definition of transverse perturbation says that $q_k(s)=e^{-s \eta_k}p_ke^{s \xi_k}$ belongs to $P$ for all 
     $s \in \RR$.  Thus 
     $$ e^{s \Ad(p_k)\xi_k}=e^{s \eta_k}q_k(s)p_k^{-1},$$
     where $s \mapsto q_k(s)p_k^{-1}$ is a smooth path in $P$ passing through $id$ when $s=0$.  Lemma \ref{lem.reparametrage} then implies  
     $$ \exp(\hat{f}_k(\hx_k).p_k^{-1}, \Ad(p_k)\xi_k)=\exp(\hat{y}_k,\eta_k).q_kp_k^{-1}.$$
     Right translation by $p_k q_k^{-1}$ gives $\hat{f}_k(\hat{z}_k).q_k^{-1}=\exp(\hat{y}_k,\eta_k)$.  
     This expression is bounded, because $(\hat{y}_k)$ is a bounded sequence, and $\eta_k$ tends to zero by definition of a 
     transverse perturbation.
\end{preuve}

 \subsection{Admissible operations}
\label{sec.notation.reduced}
 
In this section, we specialize to $X = G/P$ a parabolic model space, and define some operations on holonomy sequences specific to parabolic geometries. 
We first introduce some notation in $\lieg$.

\subsubsection{Notation in $\lieg$}
\label{sec.notation}
Let $G$ be semisimple with no compact local factors and with finite center.
We denote by $\Theta $ a Cartan involution of the semisimple Lie algebra $\lieg$.
   Associated to $\Theta$, we choose
  a Cartan subspace $\liea$, and 
  $\Phi=\{ \alpha_1,\ldots,\alpha_r\}$  a set of simple roots.  The positive
 and negative roots are denoted $\Phi^+$ and  $\Phi^-$, respectively.  
 The usual decomposition  of the Lie algebra $\lieg$ into root spaces is
 $$ \lieg=\sum \limits_{\alpha \in \Phi^-}\lieg_{\alpha} \oplus \liea \oplus \liem \oplus \sum \limits_{\alpha \in \Phi^+ } \lieg_{\alpha}.$$
 Recall that the Lie algebra $\liem$ is centralized by $\liea$, and lies in the Lie algebra $\liek$  comprising the $+1$-eigenspace of the Cartan involution $\Theta$.  
  
 We will denote by $\lien^+$ (resp.  $\lien^-$) the sum $\sum \limits_{\alpha \in \Phi^+} \lieg_{\alpha}$ 
 (resp. $\sum \limits_{\alpha \in \Phi^-} \lieg_{\alpha}$).

 The minimal parabolic subalgebra of $\lieg$ is $\liep_{min}= \liea \oplus \liem \oplus \lien^+$.
 %\sum \limits_{\alpha \in \Phi^+ } \lieg_{\alpha}$.
 A general \emph{parabolic subalgebra} $\liep$ is one containing $\liep_{min}$, and is obtained as follows (up to conjugacy in $G$): 
 there exists  $\Lambda \subsetneq \Phi$, possibly empty, such that  
      $$\liep_{\Lambda}=\sum \limits_{\alpha \in{\Lambda}^+}\lieg_{- \alpha} \oplus \liep_{min}.$$
    
   where ${\Lambda}^+$ is the set of roots in $\Phi^+$ which are in the span of $\Lambda$.   A \emph{parabolic subgroup} of $G$ is any Lie subgroup $P_\Lambda < G$ with Lie algebra $\liep_{\Lambda}$, for some $\Lambda$.  We will sometimes denote this group simply $P$ when $\Lambda$ is understood.
    
%Elements of ${\Lambda}^+$ are thus obtained
 %  as linear combinations of roots in $\Lambda$, with positive integer coefficients. 
  %\item{Generally, when $\Lambda$ is fixed once for all and there is no ambiguity, we will write $\liep$ instead of $\liep_{\Lambda}$.}   

We denote by $\lien_\Lambda^+$ the nilpotent radical of $\liep$, which equals $\sum \limits_{\alpha \in (\Lambda^+)^c} \lieg_{\alpha}$. 
Here $(\Lambda^+)^c$ stands for the positive roots written as  linear 
   combinations of roots
    in $\Phi$ involving at least one root which is not in $\Lambda$. Notice that $\lien_\Lambda^+$ is an ideal of $\lien^+$ and of
 $\liep$.  Finally, we call $\lieh_{\Lambda}$ the Lie algebra $\lieh_{\Lambda}=\liea \ltimes \lien_\Lambda^+$.

We denote by $A$, $N_{\Lambda}^+$ and $H_{\Lambda}$ the connected Lie subgroups of $G$ with Lie algebras
  $\liea$, $\lien_\Lambda^+$ and $\lieh_{\Lambda}$, respectively; they are all subgroups of $P_{\Lambda}$.
%  \item{The Levi decomposition of $\liep$ is
% $$ \liep_\Lambda \cong \lieh_\Lambda \ltimes \lien_\Lambda^+$$
% where
%  $$\lieh_{\Lambda}=\liea \oplus \liem \oplus \sum_{\alpha \in{\Lambda}^+} \left( \lieg_{\alpha} \oplus \lieg_{- \alpha}\right).$$}
%  \end{itemize}

\subsubsection{Reduced holonomy sequences}

A sequence $(p_k)$ in $P$ will be called  \emph{reduced} when it is a sequence of $H_{\Lambda}$. 

%  If $(f_k)$ is a sequence of automorphic immersions of a Cartan geometry modeled on $G/P$, we will  denote by $\overline{\Hol}(x) \subset \Hol(x)$ the  set of all  holonomy sequences at $x$ which are reduced.
 
 \begin{lemme}
\label{lemma.reducing.holonomy} 
% Let $(f_k)$ be a sequence of automorphic immersions which is equicontinuous at $x \in M$.  Then 
Any sequence $(p_k)$ in $P = P_\Lambda$ can be converted by left and right vertical perturbation to $(q_k) \in H_\Lambda$.
 \end{lemme}
 
\begin{preuve}
Consider the Levi decomposition of $P_\Lambda = S_\Lambda \ltimes N_\Lambda^+$, where $S_\Lambda$ is the connected reductive subgroup of $G$ with Lie algebra spanned by $\liea$ and the positive and negative root spaces of $\Lambda^+$. Write $p_k = s_k n_k$ according to this decomposition.  As $S_\Lambda$ is reductive, it admits a $KAK$ decomposition, according to which $s_k = l_k' a_k l_k$, with $a_k \in A = \exp(\liea)$ and $l_k, l_k' \in K$.  As $G$ has finite center, $K$ is contained in a maximal compact subgroup of $G$ and is a maximal compact subgroup of $S_\Lambda$.  Then $p_k = l_k' a_k n'_k l_k$, where $n_k' = l_k^{-1} n_k l_k \in N_\Lambda^+$.  Now $q_k = a_k n'_k$ is the desired reduced sequence.
\end{preuve}

\subsubsection{Weyl reflections}

For $X = G/P$ parabolic, these are transformations of holonomy sequences in $H_\Lambda$, which will be useful in our proof.

For any root $\alpha$, the {\it Weyl reflection} is $\rho_{\alpha}: \liea^* \to \liea^*$  with
   $$ \rho_\alpha(\xi) = \xi - \frac{2 \langle \alpha ,\xi \rangle}{\langle \alpha, \alpha \rangle} \alpha \qquad \xi \in \liea^* $$
%If $\alpha \in \Phi$, then $r_\alpha$ sends all positive roots save $\alpha$ to positive roots (REF). 
Recall that for $\alpha$ positive,  $\rho_{\alpha}$ preserves $\Phi^+ \backslash \{ \alpha \}$ and $\Phi^- \backslash \{ - \alpha \}$, assuming $2 \alpha$ is not a root (in which case, $\rho_\alpha$ preserves $\Phi^+ \backslash \{ \alpha, 2 \alpha \}$ and $\Phi^- \backslash \{ - \alpha, - 2 \alpha \}$).  Recall that whenever $\xi$ is a root, 
then $A_{\alpha \xi} = 2 \langle \alpha ,\xi \rangle/\langle \alpha, \alpha \rangle$ is an integer.  

 For any root $\alpha$, there exists $k_{\alpha} \in G$, such
 that $\Ad(k_{\alpha})$ preserves $\liea$, and the action of $\Ad(k_\alpha)$ on $\liea^*$ coincides with that of $\rho_{\alpha}$  (see \cite[Prop 6.52c]{knapp.lie.groups}).  In the sequel, 
 we will denote 
  by $r_{\alpha}$ any automorphism of $G$ such that the action induced on $\lieg$  preserves $\liea$ and 
  sends every root space $\lieg_{\beta}$ to the 
corresponding $\lieg_{\rho_{\alpha}(\beta)}$; for instance, $r_{\alpha}$ could be conjugacy by $k_{\alpha}$. 

Let  $\alpha \in \Lambda^+$.  If a root $\beta$ is a linear combination
 with integer coefficients of roots in $\Lambda$, then so is $\rho_{\alpha}(\beta)$; thus
 $\rho_{\alpha}$ preserves $\Lambda^+ \cup - \Lambda^+$.  As $\rho_\alpha$ sends all positive roots except multiples of $\alpha$ to positive roots, it also preserves $\Phi^+ \backslash \Lambda^+ = (\Lambda^+)^c$.
%% .  Then $\beta$ writes as a linear combination with positive integers
%%    of roots in $\Phi$ with one of those roots, say $\gamma$, not in $\Lambda$.  Then when writing $\rho_{\alpha}(\beta)$ as a linear combination
%%     with integers coefficients of roots in $\Phi$, the coefficient of $\gamma$ is again positive, meaning that $\rho_{\alpha}(\beta)$ is a 
%%     positive root, which moreover belongs to $(\Lambda^+)^c$.  
We conclude that for every $\alpha \in \Lambda^+$, an automorphism $r_{\alpha}$ preserves the 
connected subgroups $A$, $N_\Lambda^+$, and the identity component $P_\Lambda^0$; 
in particular, it sends sequences $(p_k)$ in $H_{\Lambda}$ to $r_{\alpha}(p_k)$ in $H_{\Lambda}$. 
Note that in general, $P_\Lambda$ may not be invariant by $r_\alpha$.
  
\subsubsection{Definition of admissible operations, perturbations}

\begin{definition} %[Admissible operations]
   \label{defi.admissibleoperation}
   Let $X=G/P$ be a parabolic variety with $P=P_{\Lambda}$.
   For $(p_k)$ a sequence of $P$, an \emph{elementary admissible operation} on $(p_k)$ is 
    of one of the three following types:
    \begin{enumerate}
     \item{A vertical perturbation of $(p_k)$.}
     \item{A transverse perturbation of $(p_k)$.}
\item{For $(p_k)$ in $H_{\Lambda}$, a Weyl reflection $r_{\alpha}$ applied to $(p_k)$,  with $\alpha \in \Lambda^+$.}
     
    \end{enumerate}
    
An \emph{admissible perturbation} of a sequence $(p_k)$ in $P$ is a sequence $(q_k)$  which is obtained from $(p_k)$ by finitely many elementary admissible operations.  
  \end{definition}
% \begin{remarque}
%   \label{rem.weyl}
  Note that the result of an admissible perturbation of a sequence $(p_k)$ of $P$ is always in $P$.  Weyl reflections are only allowed on sequences of $H_{\Lambda}$, which must be kept in mind when applying
  successive admissible operations.
%, some of them being Weyl reflections. 
  
  We conclude this section with an important remark about Weyl reflections.  We observed at the end of the last paragraph
   that a Weyl reflection $r_{\alpha}$  always coincides with the conjugacy by some element $k_{\alpha} \in G$. We also
    observed that $r_{\alpha}$  preserves the identity component $P^0$ of $P$, so that actually $k_{\alpha}$  belongs to 
     $\mathrm{Nor}_G(P^0)$, the normalizer  of $P^0$ in $G$.   
  This  normalizer  $\mathrm{Nor}_G(P^0)$ has Lie algebra 
  $\liep$ (see \cite[Lemma 3.1.3, Cor. 3.2.1(4)]{cap.slovak.book.vol1}), so that the inclusion $P \leq \mathrm{Nor}_G(P^0)$
   holds. Observe that in general, these groups need not coincide.  However, 
   when $P = \mathrm{Nor}_G(P^0)$, any Weyl reflection $r_{\alpha}(p_k)$ is 
    actually {\it a vertical perturbation of $(p_k)$}. We thus get a straightforward 
    rephrasing  of Lemma \ref{lem.stability.transverse}, namely
 
 %\begin{proposition}
  \begin{lemme}
   \label{lem.stableadmissible}
   Let $(M,\mathcal{C})$ be a parabolic geometry modeled on $X=G/P$, where $P=\mathrm{Nor}_G(P^0)$.
    Let  $x \in M$, and let $(f_k)$ be a sequence
    of automorphic immersions which is equicontinuous at $x$.
   Then if $(p_k)$ is in  $\Hol(x)$, any admissible perturbation of $(p_k)$  is in $\Hol(x)$. 
   %stable by vertical and transverse perturbations.  
  \end{lemme}
 %\end{proposition}

The case of equality $P = \mathrm{Nor}_G(P^0)$ will thus be technically more convenient, since it means that Weyl reflections on holonomy sequences again yield holonomy sequences.  It is explained in Section \ref{sec.gettingrid} why this equality may be assumed.

\section{Translation of the main theorem to the model space}
\label{sec.trans.to.model}

Via the holonomy sequences associated to an equicontinuous sequence $(f_k)$ of automorphic immersions, we can translate Theorem \ref{thm.punctualequicontinuity} to an assertion about sequences of $H_\Lambda$ acting equicontinuously with respect to segments on $X$. 

%% Actually, the proof will rely on the key Theorem  \ref{thm.induction} below, which is in itself an interesting property
%% about equicontinuity in the model space $X=G/P$.  

\begin{theoreme}
  \label{thm.induction}
  Let $X=G/P$ be a parabolic variety with $P=P_{\Lambda}$.  
  Given a sequence $(a_kn_k)$ of $H_{\Lambda}$ which, together with all of its admissible perturbations, acts equicontinuously with respect to segments on $X$,  
 the factor $(n_k)$ is bounded. 
 \end{theoreme}

Theorem \ref{thm.induction} is proved in sections \ref{sec.proof.rank1}, \ref{sec.tools} and \ref{sec.induction}.

\subsection{Derivation of Theorem \ref{thm.punctualequicontinuity} from Theorem \ref{thm.induction}}
\label{sec.derivation}
Given a sequence $(f_k)$ of automorphic immersions as in the statement of Theorem \ref{thm.punctualequicontinuity}, 
let $(p_k)$ be a holonomy sequence of $(f_k)$ at $x$. 
We can assume by Lemmas \ref{lem.stability.transverse} and \ref{lemma.reducing.holonomy} that $p_k \in H_\Lambda$ for all $k$.

 %% We saw in Lemma \ref{lem.stability.transverse} that vertical perturbations preserve ${{\mathcal Hol}}(x)$, so that applying Lemma \ref{lemma.reducing.holonomy}
 %% we can assume that $(p_k)$ is a sequence of $H_{\Lambda}$. For the sake of clearness, 

We will first deduce Theorem \ref{thm.punctualequicontinuity} \emph{under the extra assumption that $P$ equals $\mathrm{Nor}_G(P^0)$}.  
Section \ref{sec.gettingrid} explains how to dispense with this assumption.
  
Proposition \ref{prop.equivalence} ensures that $(p_k)$ acts equicontinuously with respect to segments on $X$.  
Lemma \ref{lem.stableadmissible} says that in fact every admissible perturbation of $(p_k)$ does (under our assumption $P = \mathrm{Nor}_G(P^0)$).  Now the hypotheses of Theorem \ref{thm.induction} are satisfied.  The conclusion implies that $(a_k)$ is a right vertical perturbation of $(p_k)$, which by Lemma \ref{lem.stability.transverse} also belongs to ${\mathcal Hol}(x)$.  The action of $\Ad(a_k)$ on $\lieg$ preserves the subalgebra 
  $\lien^{-}$; denote by $L_k$ the endomorphism $\left. \Ad(a_k) \right|_{\lien^-}$.
  
  \begin{lemme}
   \label{lem.bounded}
   The sequence $(L_k)$ is bounded in $\End(\lien^{-})$.
  \end{lemme}

\begin{preuve}
 The representation 
   of $\Ad(a_k)$ on $\lien^-$ is diagonalizable with eigenvalues $(\lambda_1(k),\ldots, \lambda_s(k))$. Assume for a contradiction that 
   $L_k$ is unbounded; we may assume that $\lambda_1(k)$ is unbounded, and after passing to a subsequence, that $|\lambda_1(k)| \to \infty$.  
   Taking a subsequence also allows us to assume that in $\hm$, the sequence
%there is a sequence $(\hx_k)$ in $\hm$, converging to $\hx \in \pi^{-1}(x)$, and such that 
$\hy_k=f_k(\hx_k).p_k^{-1}$ converges to $\hy$. 

For each $k$, let $\eta_k$ be in the $\lambda_1(k)$-eigenspace of $L_k$ such that $\eta_k \to \eta_{\infty} \not =0$; these can moreover be chosen in the injectivity domain of $\exp_{\hy_k}$, and such that $\eta_\infty$ is in the injectivity domain of $\exp_{\hy}$.  Set $\xi_k:=\eta_k/\lambda_1(k)$.  Because $\xi_k \to 0$, the exponential $\exp(\hx_k,\xi_k)$ is defined for sufficiently large $k$, and satisfies
    $$ f_k(\exp(\hx_k,\xi_k)).a_k^{-1}=\exp(\hy_k,\eta_k).$$
    Projecting to $M$ gives a contradiction to the equicontinuity of $(f_k)$ at $x$: $\pi(\exp(\hx_k,\xi_k)) \to x$, while
     $\pi(\exp(\hy_k,\eta_k)) \to \pi(\exp(\hy,\eta_{\infty})) \not = \pi(\hy)$.
\end{preuve}

% We stick to the notations introduced in the proof above.  It follows from Lemma \ref{lem.bounded} that considering  
Now again passing to a
 subsequence of $(f_k)$, we may assume that  $L_k$ tends to some $L \in \mbox{End}(\lien^{-})$.
   Let $K \subset \hm$ be a compact set containing both sequences $(\hx_k)$ and $(\hy_k)$, and let ${\mathcal U}$ and ${\mathcal V}$ be relatively compact neighborhoods of $0$ in $\lien^-$, such that:
  \begin{enumerate}
  \item{$L_k(\overline{\mathcal U}) \subset {\mathcal V}$ for every $k \in \NN$.}
\item{For every $\hz \in K$, the map $\Phi_{\hz} : u \mapsto \pi(\exp(\hz,u))$ 
is defined on  $\overline{\mathcal U}$ and 
$\overline{\mathcal V}$,
 and is a diffeomorphism from ${\mathcal U}$ and ${\mathcal V}$ onto their respective images. 
}
  \end{enumerate}
 There exists an open neighborhood $U$ of $x$, such that  $U \subseteq \Phi_{\hz}({\mathcal U})$ for $\hz \in K$ close enough to $\hx$.  Then define the smooth map $h : U \to M$
  by $h=\Phi_{\hy} \circ L \circ \Phi_{\hx}^{-1}$.  Because $L_k$ converges smoothly to
   $L$, and since on $U$, for $k$ large enough,
   $$ f_k=\Phi_{\hy_k}\circ L_k \circ \Phi_{\hx_k}^{-1},$$
  $(f_k)$ converges smoothly to $h$ on $U$. Thus Theorem \ref{thm.punctualequicontinuity} is proved.

\subsection{Justification of the assumption $P=\mathrm{Nor}_G(P^0)$} 
\label{sec.gettingrid}

Let $(f_k)$ be a sequence of automorphic immersions as in Theorem \ref{thm.punctualequicontinuity}.  In general 
$P \leq \mathrm{Nor}_G(P^0)$, and they have the same Lie algebra, as remarked above (again, see \cite[Lemma 3.1.3, Cor. 3.2.1(4)]{cap.slovak.book.vol1}).  Thus $P' = \mathrm{Nor}_G(P^0)$ is an isogenous supergroup of $P$.  The following lemma gives a general procedure for inducing a Cartan geometry modeled on $G/P$ to one modeled on $G/P'$, with 
respect to which the automorphism group behaves nicely.

 \begin{lemme}
 \label{lem.adjust}
 Let ${\mathcal C} = (\widehat{M},\omega)$ be a Cartan geometry on the manifold $M$, modeled on $X=G/P$. 
 Let $P^{\prime} < G$ be a closed subgroup, with $P \leq P'$ and $(P')^0 = P^0$. 
 Then there exists a Cartan geometry ${\mathcal C}^{\prime}=({\widehat{M}}^{\prime},\omega^{\prime})$ on the manifold $M$, modeled on
 $X^{\prime}=G/P^{\prime}$, such that:
 \begin{enumerate}
  \item Every automorphic immersion of $(M,{\mathcal C})$ is an automorphic immersion of $(M,{\mathcal C}')$.
 \item The corresponding inclusion of $\Aut(M,{\mathcal C})$ into $\Aut(M,{\mathcal C}^{\prime})$ is a homeomorphism onto a closed subgroup
 with respect to the Lie topologies on each.
 %; similarly, the automorphic immersions of $(M,\mathcal{C})$ defined on an open $U \subseteq M$ embed homeomorphically as a closed subset of the automorphic immersions of $(M,\mathcal{C}')$ defined on $U$, in the respective Lie topologies.   
 
% \item The topology induced on $\Aut(M,{\mathcal C})$ by the Lie topology of $\Aut(M,{\mathcal C}^{\prime})$  coincides with the Lie topology of $\Aut(M,{\mathcal C})$.
 \end{enumerate}

  %Moreover, 
   %the Lie topology of $\Aut(M,{\mathcal C}^{\prime})$ coincides with the Lie topology of $\Aut(M,{\mathcal C})$.
   %of $\Aut(M,{\mathcal C}^{\prime})$ 
%   Moreover, if $\Aut(M,{\mathcal C})$ and $\Aut(M,{\mathcal C}^{\prime})$ are endowed
%    with their respective Lie topologies, the inclusion map $\iota: \Aut(M,{\mathcal C})$
%  % $P \subset P^{\prime} \subset G$ with $dim(P)=dim(P^{\prime})$.  
  \end{lemme}

\begin{preuve}
% We consider the product space ${\mathcal M}={\widehat{M}} \times P^{\prime}$.  There are two natural and pairwise commuting actions on ${\mathcal M}$.
%  A left-action of $P$ given by
% and a right action of $P^{\prime}$ defined by $(\hx,q).q_0=(\hx,qq_0)$.  The action of  $P$ on ${\mathcal M}$
%  is free and proper, so that we can consider 
The bundle $\widehat{M}'$ is obtained as the quotient ${\widehat{M}} \times_P P^{\prime}$, where $P$ acts diagonally by $p.(\hx,q)=(\hx.p^{-1},pq)$, freely and properly.  There is an obvious commuting right $P'$-action on $\mathcal{M} = \widehat{M} \times P'$, which descends to $\widehat{M}'$, making it a $P^{\prime}$-principal bundle over $M$.

To construct the Cartan connection on $\widehat{M}'$, we first build a one-form $\tilde{\omega} \in \Omega^1(\mathcal{M},\lieg)$.  For $(\xi,u) \in T_{(\hx,q)}\mathcal{M}$, let
     $$ \tilde{\omega}_{(\hx,q)}(\xi,u):= \Ad(q^{-1})\omega_{\hx}(\xi)+ (\omega_{P^{\prime}})_q(u).$$
     where $\omega_{P^{\prime}}$ is the Maurer-Cartan form of $P^{\prime}$.  
It is readily checked that $\tilde{\omega}$ satisfies the equivariance relation 
$(R_{p})^* \tilde{\omega}= \Ad(p^{-1}) \circ \tilde{\omega}$ for every $p \in P^{\prime}$, and that it is invariant under the diagonal action of 
$P$ on $\mathcal{M}$.
Moreover
$$\tilde{\omega}_{(\hx,q)}(T_{\hx}{\widehat{M}} \times \{0\})= \Ad(q^{-1}) \circ \omega_{\hx}(T_{\hx}{\widehat{M}})=\lieg$$
 showing  that
$\tilde{\omega}: T{\mathcal M} \to \lieg$ is onto at each point.  

For $X \in \liep$, let $X^\ddag \in \mathcal{X}({\widehat{M}})$
% induced by the right action of $e^{tX}$ on ${\widehat{M}}$, 
be as in Definition \ref{defi.cartan.geometry}, and let $\gamma$ be the curve 
$$\gamma(t)=e^{tX}.(\hx,q)=(\hx.e^{-tX},e^{tX}q).$$  
Then 
$$\tilde{\omega}(\gamma^{\prime}(t))= \Ad(q^{-1})\circ \omega_{\hx}(-X^\ddag)+ \Ad(q^{-1})X=0$$ 
since
 $\omega(X^\ddag) \equiv X$.  Hence the kernel of $\tilde{\omega}_{(\hx,q)}$ contains the tangent space to the $P$-orbits on ${\mathcal M}$; by a dimension argument, these spaces are equal.  We infer that $\tilde{\omega}$ induces a $1$-form $\omega^{\prime} \in \Omega^1({\widehat{M}}^{\prime},\lieg)$, which is the desired Cartan connection on ${\widehat{M}}^{\prime}$. 
  
% It remains to prove the inclusion $Aut(M,{\widehat{M}},\omega) \subset Aut(M,{\widehat{M}}^{\prime},\omega^{\prime})$.  
We prove point (1) for $f \in \Aut(M,{\mathcal C})$. The argument for automorphic immersions is similar.
Let $\hat{f}$ be the lift of $f$ to $\widehat{M}$, and define 
 ${\tilde{f}} : {\mathcal M} \to {\mathcal M}$ by ${\tilde{f}}(\hx,q)=({\hat f}(\hx),q)$.  
 The $P$-equivariance of ${\hat f}$ gives the equivariance relation 
   $p.{\tilde{f}}(\hx,q)={\tilde{f}}(p.(\hx,q))$;  obviously, 
 ${\tilde{f}}((\hx,q).p^{\prime})={\tilde{f}}(\hx,q).p^{\prime}$ for every $p^{\prime} \in P^{\prime}$.  Thus
    ${\tilde{f}}$ induces a bundle morphism ${\hat f}^{\prime}$ of ${\widehat{M}}^{\prime}$ covering $f$. 
    
To prove that $f \in \Aut(M,\mathcal{C}')$,
it remains to check 
%to get the inclusion  $f \in Aut(M,{\widehat{M}}^{\prime}, \omega^{\prime})$,  is 
that ${\hat f}^{\prime}$ preserves 
    $\omega^{\prime}$. 
    To this end, we compute ${\tilde{f}}^* \tilde{\omega}$ and show that it coincides with $\tilde{\omega}$:  
    $$ \tilde{\omega}_{({\hat f}(\hx),q)}(D_{\hx}{\hat f}(\xi),u)= \Ad(q^{-1}) \circ \omega_{{\hat f}(\hx)}(D_{\hx}{\hat f}(\xi))+ (\omega_{P^{\prime}})_q(u)$$
    but $\omega_{{\hat f}(\hx)}(D_{\hx}{\hat f}(\xi))=\omega_{\hx}(\xi)$ because $f \in \Aut(M,\mathcal{C})$.  Finally,
    $$ \tilde{\omega}_{({\hat f}(\hx),q)}(D_{\hx}{\hat f}(\xi),u)= \Ad(q^{-1})\omega_{\hx}(\xi) + (\omega_{P^{\prime}})_q(u)=
    \tilde{\omega}_{(\hx,q)}(\xi,u)$$
    as desired, so (1) is proved.
    
%% Before establishing the topological properties of the inclusion $\Aut(M,{\mathcal C}) \leq \Aut(M,{\mathcal C}^{\prime})$, let us make
%% a few elementary observations.  

There is a natural $P$-equivariant, proper embedding $j : (\hm, \omega) \to (\hm',\omega')$ defined by $j(\hx):=[(\hx,e)]$, the $P$-orbit in $\mathcal{M}$ of $(\hx,e)$.  
%Here, we denote by $[(\hx,e)]$ the class in $\hm'$ of the element $(\hx,e) \in {\mathcal M}$.  It is straigthforward to check that $j$ is $P$-equivariant, and that $j^* \omega'=\omega$. Moreover,  
For $f \in \Aut(M, {\mathcal C})$ with respective lifts ${\hat f}$ and ${\hat f}'$ to $\hm$ and $\hm^{\prime}$, we have $j \circ {\hat f} = {\hat f}' \circ j$.

Now consider a sequence $f_k \in \Aut(M,{\mathcal C})$ converging for the Lie topology of $\Aut(M,{\mathcal C}')$ to an automorphism $f$. 
%By the characterization of the Lie topology given in 
By Kobayashi's theorem (Thm \ref{theo.kobayashi}), the sequence of lifts 
${\hat f}_k'$ 
%constructed at the begining of the proof 
converges in the $C^{\infty}$-topology of $\hm'$ to a diffeomorphism ${\hat f}'$, which clearly preserves $\omega'$.  Properness of $j$ implies that $j(\widehat{M})$ is closed.  Then ${\hat f}'$ preserves $j(\hm)$, because every ${\hat f}_k$ does.  Thus ${\hat f}_k=j^{-1} \circ {\hat f}_k' \circ j$  converges smoothly on $\hm$
    to ${\hat f}:=j^{-1} \circ {\hat f}' \circ j$, which preserves $\omega$ and
     covers $f$.  It follows that $f \in \Aut(M, {\mathcal C})$, and by Kobayashi's theorem, $f_k \to f$ in the Lie topology of $\Aut(M,{\mathcal C})$.  We conclude moreover that $\Aut(M,\mathcal{C})$ is closed in the Lie topology of $\Aut(M,\mathcal{C}')$.  
 
 Conversely, given $f_k \rightarrow f$ in the Lie topology of $\Aut(M,{\mathcal C})$, with $f \in \Aut(M,\mathcal{C})$, the lifts ${\hat f}_k \to {\hat f}$ smoothly on $\hm$. These correspond, as in the proof of (1), to automorphisms ${\hat f}_k'$ and ${\hat f}'$ of $(\hm',\omega')$
%, and the smooth convergence  $j \circ {\hat f}_k \to j \circ {\hat f}$  means that 
with ${\hat f}_k' \to {\hat f}'$ on $j(\hm)$.  For any $\hy \in \hm'$, there exists $p' \in P'$ such that $\hy.p' \in j(\hm)$.
       It follows by Theorem \ref{theo.kobayashi} (3) that ${\hat f}_k' \to {\hat f}'$ smoothly on each connected component of $\hm'$; in other words, $f_k \to f$ holds in the Lie topology of $\Aut(M,{\mathcal C}')$.  Thus $\Aut(M,\mathcal{C}) \hookrightarrow \Aut(M,\mathcal{C}')$ is a homeomorphism onto its image with respect to the Lie topologies on each group.
\end{preuve}
%  Denoting by $\theta$ the projection from $ {\mathcal M}$ to ${\widehat{M}}^{\prime}$, the formula $\hx \mapsto \theta(\hx,e)$ yields 
%  a $P$-equivariant map $j : {\widehat{M}} \to {\widehat{M}}^{\prime}$ which is a diffeomorphism on its image.  We will denote in the following ${\widehat{M}}_0=j({\widehat{M}})$.
%  
     
Now, given a sequence $(f_k)$ as in Theorem \ref{thm.punctualequicontinuity}, Lemma \ref{lem.adjust} with $P^{\prime}=\mbox{Nor}(P^0)$ allows us to consider
  $(f_k)$ as a sequence of automorphic immersions of $(M,{\mathcal C}')$, modeled on $X^{\prime}=G/P^{\prime}$.  The proof of Section \ref{sec.derivation} says that $(f_k)$ converges smoothly on $M$ to a smooth map $h$.  We have thus shown that Theorem  \ref{thm.induction} implies Theorem \ref{thm.punctualequicontinuity}.  

\subsection{Derivation of  Theorem \ref{thm.lietransfogroup}}
\label{sec.proofcorollary}
Let $f_k \in \Aut(M,{\mathcal C})$ converge to $h \in \operatorname{Homeo}(M)$ in the $C^0$ topology.  The aim is to show that $h \in \Aut(M,{\mathcal C})$, and $f_k \rightarrow h$ in the Lie topology on $\Aut(M,{\mathcal C})$.  
 
 By Lemma \ref{lem.adjust} point (2), we may assume that the model space $G/P$ satisfies $P=\mathrm{Nor}_G(P^0)$.  % We can reproduce the proof  of Theorem \ref{thm.punctualequicontinuity} given 
As in Section \ref{sec.derivation}, $(f_k)$ admits a holonomy sequence $a_k \in A$ at any $x \in M$, such that $L_k= \left. \Ad(a_k) \right|_{\lien^-}$ is bounded in $\End(\lien^-)$.  Moreover, in the notation of Section \ref{sec.derivation},  there is a neighborhood $U$ of $x$ such that
   for any accumulation point $L$ of $(L_k)$ in $\End(\lien^-)$, a subsequence of $(f_k)$ converges to $\Phi_{\hy} \circ L \circ \Phi_{\hx}^{-1}$ on $U$.  Then $\left. L \right|_U =\Phi_{\hy}^{-1} \circ h \circ \Phi_{\hx}$, so $L_k \rightarrow L$.  Because $h$ is a homeomorphism, $L$ is injective around $0$, hence $L \in \GL(\lien^-)$.  As a consequence, $(a_k)$ converges in $P$.
      
Now we have $\hat{f}_k(\hx_k).a_k^{-1}=\hy_k \rightarrow \hat{y}$ with $(a_k)$ also converging, so $f_k(\hx_k)$ tends to some point $\hz$. 
As $\hat{x}_k \rightarrow \hat{x}$, for sufficiently large $k$, $\hx=\exp(\hx_k,\xi_k)$, with $\xi_k \to 0$ in $\lieg$.  Now $\hat{f}_k(\hx)=\exp(\hat{f}_k(\hx_k),\xi_k)$, so $f_k(\hx) \to \hz$.  By Theorem \ref{theo.kobayashi} (3), ${\hat f}_k $ and the inverses ${\hat f}_k^{-1}$ both converge $C^{\infty}$ on $\hm$ to smooth maps $\hat{f}$ and ${\hat g}$, which obviously satisfy $\hat{f} \circ \hat{g} =id$.
 It is easy to see that $\hat{f}$ is a bundle automorphism of $\hm$ preserving $\omega$.  
It lifts $h$, hence $h \in \Aut(M,{\mathcal C})$.
       Finally, because $\hat{f}_k \rightarrow \hat{f}$ smoothly on $\hm$, Theorem \ref{theo.kobayashi} (2) gives that $f_k \to h$ in the Lie topology on $\Aut(M,{\mathcal C})$.

\section{Proof of Theorem \ref{thm.induction} in rank one}
\label{sec.proof.rank1}
 
Our proof of Theorem \ref{thm.induction} will proceed by induction on $\rk_\RR(G)$. 
The essential arguments for the base case, $\rk_\RR(G) =1$, are in the paper
\cite{frances.lfrank1} by the first author.  For the convenience of the reader,
the proof is presented here in a manner consistent with our terminology and notation.
 Theorem \ref{thm.induction} in this rank one case will actually be a consequence of the following proposition.

%% By Lemma \ref{lemma.reducing.holonomy} and Proposition \ref{prop.equivalence}, it suffices in order to prove Theorem \ref{thm.equicontinuity} in rank one to prove the following proposition:
 
  \begin{proposition}
  \label{prop.nonequicontinuous.rk1}
  Let $X=G/P$ be a parabolic space, with $\rk_\RR(G) =1$.  If $p_k=a_kn_k$ is a sequence of $A \ltimes N^+$ such that $(n_k)$ is unbounded, then $(p_k)$
    does not act equicontinuously with respect to segments.
  \end{proposition}
  
Recall the notation of Section \ref{sec.notation}.  The rank one Lie algebra can be decomposed as a vector space direct sum of subalgebras 
$\lieg=\lien^{-}\oplus \liea \oplus \liem \oplus \lien^+$.  The Lie algebra $\lien^-$ (resp. $\lien^+$) is abelian
   if $\lieg=\mathfrak{o}(1,n)$, and nilpotent of index $2$, with center  of respective dimension $1$, $3$ and $7$ if
    $\lieg$ is $\mathfrak{su}(1,n)$, $\mathfrak{sp}(1,n)$ or $\mathfrak{f}_4^{-20}$. In all cases, ${\mathfrak z}^-$ (resp.  ${\mathfrak z}^+$)
    will
     denote the center of $\lien^-$ (resp. $\lien^+$). The nonequicontinuity will be observed on a restricted class of segments, namely
     those $[\xi]$ with
     $$ \xi \in {\bf Q} =\{ \Ad(p) u \ | \ u \in {\mathfrak z}^-, \ p \in P   \}.$$
     This set of segments will be denoted $[{\bf Q}]$ and corresponds to  conformal circles when $\lieg=\mathfrak{o}(1,n)$, and to chains and 
     their generalizations in the other rank one 
     models.  We will  adopt the 
     notation $\dot{\bf Q}$ (resp.  $[\dot{\bf Q}]$) for ${\bf Q} \setminus \{ 0 \}$ (resp. $[{\bf Q}] \setminus \{ [o] \}$).

We now recall two results from \cite{frances.lfrank1} regarding these distinguished segments.
  
  \begin{lemme}[\cite{frances.lfrank1}, Lemme 2]
   Let $([\alpha_k])$ be a sequence of segments in $[{\bf Q}]$.  If $[\alpha_k]$ tends to $[o]$ for the Hausdorff topology, then $L([\alpha_k]) \to 0$.
  \end{lemme}

  \begin{lemme}[\cite{frances.lfrank1}, Proposition 1, (ii)]
   There exists a continuous section $s: [\dot{\bf Q}] \to \dot{\bf Q}$.  In other words, if a sequence of segments $([\alpha_k])$ tends to
    a segment $[\beta] \not = [o]$, there is a convergent sequence $(\xi_k)$  in $\lieg$  such that $[\alpha_k]=[\xi_k]$.
  \end{lemme}

  By these two lemmas, if we can find a sequence of segments $[\alpha_k]$ in $[\dot{\bf Q}]$ tending to $[o]$, such that
    $p_k.[\alpha_k]$ tends to $[\beta] \in [\dot{\bf Q}]$ (maybe considering a subsequence), then $(p_k)$ does not act equicontinuously
  with respect to segments.
  
  The group $A$ has exactly two fixed points on $X$, namely $o$ and another point $\nu$. To better understand the action of $P$ on $[{\bf Q}]$, it is  convenient to work in the chart $\rho: \lien^+ \to X \setminus \{  o\}$, 
  given by $\rho(x)=e^{x}.\nu$.  In this chart, elements of $P$ act as affine transformations, and  segments $[\alpha] \in [\dot{\bf Q}]$ 
   coincide with  half-lines $[x,u)=\{  x+tu\ | \ t \in \RR \}$, where $x \in \lien^{+}$ and
    $u$ is a unit vector in ${\mathfrak z}^+$  (for any given norm in $\lieg$ which is invariant by the Cartan involution).  More precisely, the action of $A$ in the chart $\rho$ is linear, and is equivalent to the adjoint action on $\lien^+$, and 
   the action of  an element $n=e^{\xi}$, $\xi \in \lien^+$,  is given, by the Baker-Campbell-Hausdorff formula, by  $x \mapsto (\mbox{Id} + (\ad \xi)/2)(x)+\xi$, $\forall x \in \lien^+$.
   
   Now, let us write $n_k=e^{v_k}$.  By assumption, $(v_k)$ is an unbounded sequence in $\lien^+$.
 We claim there is an unbounded sequence $(x_k)$ in $\lien^+$ such that 
  \begin{equation}
  \label{eq.uk}
  x_k+\frac{1}{2} [v_k,x_k]+v_k=0
  \end{equation}
%  We claim that $(x_k)$ is unbounded.  
  To see this, decompose $\lien^{+}$ as a direct sum $\lien^+=\lieh \oplus \mathfrak{z}^+$  
  (observe that $\lieh=\{0\}$ when $\lieg=\mathfrak{o}(1,n)$).
  %Denoting by $\overline{x}_k$, respectively $\overline{v}_k$, and $\tilde{x}_k$, respectively $\tilde{v}_k$, the components of $x_k$, respectively $v_k$, on $\lieh$ and $\mathfrak{z}^+$,
  Split Equation (\ref{eq.uk}) into two equations in $\lieh$ and $\mathfrak{z}$, namely 
  $$ \overline{x}_k+\overline{v}_k=0 $$
  where $\overline{x}_k$ and $\overline{v}_k$ are the components of $x_k$ and $v_k$ on $\lieh$, respectively, and
 $$ \tilde{x}_k+\frac{1}{2}[\overline{v}_k,\overline{x}_k]+\tilde{v}_k=0$$
  where $\tilde{x}_k$ and $\tilde{v}_k$ are the components of $x_k$ and $v_k$ on $\mathfrak{z}^+$.
  If $(\overline{v}_k)$ is unbounded, then so is $(\overline{x}_k)$, and the same is true for $(x_k)$.  If $(\overline{v}_k)$ is bounded, then $(\tilde{v}_k)$ is unbounded because $(v_k)$ is unbounded.  
  This forces $(\tilde{x}_k)$ to be unbounded.
% because if it were not, the sum   $ \tilde{x}_k+\frac{1}{2}[\overline{v}_k,\overline{x}_k]$ would be bounded. But this sum is equal to $-\tilde{v}_k$: contradiction.
  
  We can now conclude the proof of Proposition \ref{prop.nonequicontinuous.rk1}. 
  Since $a_kn_k(x_k)=0$, then for $\xi$  of norm $1$   in $ {\mathfrak z}^+$, the sequence of half-lines $[x_k,\xi)$ 
  is mapped to $[0,\xi)$ by $(p_k)$.
    Now, after taking a subsequence, $x_k/|x_k|$ tends to $\xi_{\infty}$.  Thus for $\xi \not = -\xi_{\infty}$, 
     the sequence of half-lines $[x_k,\xi)$  goes to infinity in the chart $\rho$, which means that the corresponding sequence of segments 
     $[\alpha_k]$ tends to $[o]$ in $X$. On the other hand, $p_k([\alpha_k])$ is equal to a constant segment $[\alpha] \not = [o]$, and the
     non equicontinuity of $(p_k)$
     with respect to segments follows.

 \section{Tools for the induction step: sliding along root spaces}
\label{sec.tools}

The proof in the previous section for $\rk_\RR(G) =1$ relies heavily on the fact that the action of $P$ on the complement of its fixed
point $o \in G/P$ is by affine transformations.  In higher rank, the $P$-action on $G/P$ is a compactification of
 an affine action, but no longer a one point compactification.
%, and the situation is much more subtle to understand.  This is why a direct
 % proof of Theorem \ref{thm.induction} for all ranks seems hard to obtain, and justify our choice of a proof 
This difference creates significantly more complexity, which motivates our choice to prove Theorem \ref{thm.induction} by induction rather than directly in arbitrary rank.  

  %%  \subsection{ Admissible vertical perturbations}
  %%  An admissible vertical perturbation of a sequence $p_k=a_kn_k$ in $H_{\Lambda}$, is obtained
%%  by left-multiplying by a bounded sequence of $A$, and
  %% right-multiplying by a bounded sequence of $N_{\Lambda}^+$. Admissible vertical perturbations are obviously particular instances of vertical perturbations, so they preserve $\Hol(x)$.

The tools developed in this section build on those of Sections \ref{sec.vertical.transverse} and \ref{sec.notation.reduced}, with the purpose of simplifying holonomy sequences. 

\subsection{Essential range of $(p_k)$}

The group exponential of $G$ restricts to a diffeomorphism of $\liea$ onto $A$ by definition.  Moreover, $\Ad N_\Lambda^+$ is unipotent, and $Z(G) \cap N_\Lambda^+ = 1$, so $N_\Lambda^+$ is simply connected; thus $\exp$ restricts to a diffeomorphism $\lien_\Lambda^+ \rightarrow N_\Lambda^+$.  

Fix an ordering $\alpha_1 > \cdots > \alpha_r$ of $\Phi$, and endow $\Phi^+$ with the lexicographical ordering.  Then we obtain exponential coordinates $\ln a = (Z^1, \ldots, Z^r)$ on $A$ and $\ln n = Y = (Y^\alpha)_{\alpha \in (\Lambda^+)^c}$, where $Y^\alpha$ is a vector in $\lieg_\alpha$, on $N_\Lambda^+$.

\begin{proposition}
\label{prop.trivial.or.unbounded}
Let $p_k = a_k n_k \in H_\Lambda$ with exponential coordinates $((Z^i_k), (Y^\alpha_k))$.  Then up to vertical perturbation of $(p_k)$, we may assume each component sequence $(Y^\alpha_k)$ is either trivial or unbounded.  
\end{proposition}

\begin{preuve}
%% Suppose $(Z_k^i)$ is bounded.  Then left multiplication of $p_k = a_k n_k$ by $e^{- Z_k^i}$ results in $p_k' = a_k' n_k$ with new coordinate $((Z^i)'_k)$ bounded.  All other coordinates are unchanged because $A$ is abelian.
The group $N_\Lambda^+$ is nilpotent; write the lower central series 
$$N_\Lambda^+ = N^{(0)} \rhd N^{(1)} \rhd \cdots \rhd N^{(d)} \rhd id $$
Each $\lien^{(i)}/\lien^{(i+1)}$ is abelian and can be spanned by a direct sum of certain root spaces; denote the corresponding set of roots by $\Sigma^{(i)}$.  Let $\Pi \subset (\Lambda^+)^c$ be the set of roots $\alpha$ with $(Y^\alpha_k)$ bounded.  We first multiply $p_k$ on the right by $e^{- Y^\alpha_k}$ for all $\alpha \in \Pi \cap \Sigma^{(0)}$, in any order.  The Baker-Campbell-Hausdorff formula implies that the resulting exponential coordinates $((Y')^\alpha_k)$ are trivial or bounded for all $\alpha \in \Pi \cap \Sigma^{(0)}$.  Then proceed sequentially through $\Pi \cap \Sigma^{(i)}$ for $i=1, \ldots d$ to obtain $(p_k')$ satisfying the conclusion of the proposition. 
%, that each $((Y^\alpha)'_k)$ is trivial or unbounded. 
 \end{preuve}

We remark that $(Z_k^i)$ can also be assumed trivial or bounded by a similar argument, which is not given because this fact is not needed below.  

\begin{definition}
Let $p_k = a_k n_k \in H_\Lambda$ with exponential coordinates $((Z^i_k), (Y^\alpha_k))$. 
The \emph{essential range} of $(p_k)$, denoted $ER(p_k)$, is the set of roots $\lambda \in (\Lambda^+)^c$ for which the component sequence $(Y^\lambda_k)$ is unbounded.   
\end{definition}

\subsection{Transverse and vertical sliding along root spaces}

In our proof by induction on the rank of $G$, the goal will be, given a sequence $(p_k)$ in $H_\Lambda$, to obtain roots in the essential range of $(p_k)$ that belong to a lower-rank subspace of the span of $\Phi$. 
More precisely, given $\lambda \in ER(p_k)$ such that $\lambda$ has nontrivial component on some $\alpha \in (\Lambda^+)^c$, we will perform admissible perturbations on $(p_k)$ to obtain a new sequence $(q_k) \subset H_\Lambda$ with $\lambda - \alpha \in ER(q_k)$.  Such a manipulation is possible only under some circumstances, which are enunciated in Propositions \ref{prop.transverse.sliding} and \ref{prop.vertical.sliding} below.  First, the following proposition holds the basic Lie-algebraic calculations that make our ``sliding along $\lieg_{-\alpha}$'' procedure work.

\begin{proposition}
\label{prop.minus.alpha.conjugation}
Assume that $\alpha, \nu, \nu + \alpha \in \Phi^+$.  Given a sequence $(Y_k)$ in $\lien^+$ with $(Y_k^{\nu + \alpha})$ unbounded, there exists $\xi_k \rightarrow 0$ in $\lieg_{- \alpha}$ such that
\begin{enumerate}
\item $[\xi_k,Y_k^{\nu + \alpha}] = [\xi_k,Y_k]^\nu$ is unbounded
\item  $\left( \Ad (e^{\xi_k}) Y_k \right)^\nu$ is unbounded
  \end{enumerate}
\end{proposition}

\begin{preuve}
The bilinear map $\lieg_{- \alpha} \times \lieg_{\nu+\alpha} \to \lieg_{\nu}$ induced by the bracket is nondegenerate; we recall the proof of 
this fact for real semisimple Lie algebras.  
Denote $B$ the Killing
 form on $\lieg$; $\Theta$ the  Cartan involution as in section \ref{sec.notation}; and $H_{\nu + \alpha} \in \liea$ the dual with respect to $B$ of $\nu + \alpha$.  
%such that $B(H_{\alpha},H)=\alpha(H)$ for all $H \in \liea$. 
 Then, given $Y \in \lieg_{\nu + \alpha}$ nonzero, $[\Theta(Y),Y]=B(\Theta(Y),Y)H_{\nu + \alpha}$. Rescaling $Y$ if necessary, the vectors  
$Y$, $\Theta(Y)$ and $[\Theta(Y),Y]=H$  form an $\mathfrak{sl}_2$-triple.  Consider $V=\oplus_{k \in {\BZ}} \lieg_{- \alpha + k(\nu + \alpha)}$, which
 is an $\mathfrak{sl}_2$-module.  If $[\lieg_{-\alpha},Y]$ were zero, then 
 $V^{\prime}=\oplus_{k \leq 0} \lieg_{- \alpha + k(\nu + \alpha)}$ would be a submodule with highest weight  $-\alpha(H)$, which implies $\alpha(H) < 0$.  
%In particular, since $\alpha(H)<0$,  we woud have  $0<\nu(H)$.
  On the other hand, $V/V^{\prime}$ is also an $\mathfrak{sl}_2$-module with lowest weight $\nu(H) = - \alpha(H) + (\nu + \alpha)(H) > 0$, which is impossible.

  Given $Y \in \lieg_{\nu+\alpha}$, $|Y|=1$ (for any norm on $\lieg$), let
  $$ m(Y)=\max_{X\in \lieg_{-\alpha}, \, |X|=1}|[X,Y]| >0.$$  
Then $\inf_{Y \in \lieg_{\nu+\alpha}, \, |Y|=1}m(Y) \geq c>0$.  In particular, there exist
  $\xi_k \in \lieg_{-\alpha}$,  $|\xi_k|=1$ such that
  $$| [\xi_k,Y_k]^\nu| = |[\xi_k,Y_k^{\nu + \alpha}]|=m \left( \frac{Y_k^{\nu+\alpha}}{|Y_k^{\nu+\alpha}|} \right) |Y_k^{\nu+\alpha}| \geq c |Y_k^{\nu+\alpha}|$$
  is unbounded.  Observe that replacing $\xi_k$ by ${\xi_k}/{|Y_k^{\nu+\alpha}|^{1/2}}$ gives the same conclusion with the extra property that $\xi_k \to 0$.  Now (1) is proved.

  The conjugates in (2) are given, for some $m \in \NN$, by 
    $$ \Ad(e^{\xi_k}) Y_k = Y_k^{\prime}=\sum_{j=0}^{m} \frac{1}{j!}(\ad \ \xi_k)^j(Y_k)$$

After replacing $\xi_k$ with $s \xi_k$, the $\nu$ components are
 $$Y_k^{' \nu} =  \sum_{j=0}^m \frac{s^j}{ j!} (\ad \xi_k)^j(Y_k^{\nu + j \alpha}) $$
 
    From (1), the $\nu$ components of the terms corresponding to $j=1$ form an unbounded sequence.  The following lemma shows that replacing 
    $\xi_k$ by $s \xi_k$, with a suitable $s \in (0,1]$, makes the components $(Y_k^{' \nu})$
    unbounded too. \end{preuve}

  \begin{lemme}
  Let $(u_0(k)),\ldots, (u_m(k))$  be $m$ sequences in a finite dimensional vector space $V$. Assume that one of the sequences 
  $(u_j(k))$ is unbounded. Then for a suitable
   choice of $s \in (0,1]$, the sequence $u_0(k)+s u_1(k)+s^2 u_2(k)+\cdots+s^m u_m(k)$ is unbounded.
 \end{lemme}
\begin{preuve}
There exist $(m+1)$ values of $s$ in $(0,1]$, say $s_0,\ldots,s_m$, such that the vectors $v_i=(1,s_i,\ldots,s_i^m)$ 
form a basis of $\RR^{m+1}$.  Let $|\cdot|$ be any norm on $V$.  Then on the vector space of linear maps $\mathcal{L}(\RR^{m+1},V)$, we have 
 two norms: 
 $$ ||f||_1=\sup_{|v|=1}|f(v)|$$
 and
 $$ ||f||_2=\max_{i=0,\ldots,m}|f(v_i)|.$$
 If $f_k$ denotes the linear map $(\lambda_0,\ldots,\lambda_m) \mapsto \lambda_0 u_0(k)+\cdots+\lambda_m u_m(k)$, then 
  $(||f_k||_1)_{k \in \NN}$ is unbounded (which is the case under the hypothesis of the lemma) if and only if  $(||f_k||_2)_{k \in \NN}$ is unbounded.  The lemma follows.
\end{preuve}

Define $\Phi^+_{max} \subset \Phi^+$ to be the subset comprising the positive roots in which all $\alpha_i \in \Phi$ 
occur with a positive coefficient.  Observe  that this set is nonempty only when $G$ is simple.

\begin{proposition}(Transverse sliding)
\label{prop.transverse.sliding}
Let $p_k=a_k n_k \in H_\Lambda$ with $ER(p_k) \subseteq \Phi^+_{max}$, and assume $(p_k)$ and all its admissible perturbations act equicontinuously with respect to segments on $G/P$.
 Let $\alpha \in (\Lambda^+)^c$ such that 
% \begin{itemize}
%\item 
for all $\lambda \in ER(p_k)$, for all $l \geq 0$, if $\lambda - l \alpha$ is a root, then it belongs to $(\Lambda^+)^c$.
%\item 
%\end{itemize}
Suppose $\alpha+\nu \in ER(p_k)$ for some $\nu \in \Phi^+$.  Then vertical and transverse perturbation of $(p_k)$ yields $q_k=a_kn_k' \in H_\Lambda$ such that $\nu \in ER(q_k)$.
  \end{proposition}

% 
% By our assumption that $\alpha_{j_0}$ does not enter i the decomposition of $\alpha$, we get that $\alpha \not \in \Phi_{max}^+$, hence 
%   the component  $Y_k^{\alpha}$ is trivial.
  %$Y_k^{\alpha}$ of $$
\begin{preuve}
If $(Y_k^\nu)$ is unbounded, there is nothing to do.  By Proposition \ref{prop.trivial.or.unbounded}, we may assume after a vertical perturbation that $Y_k^\lambda$ is trivial for all $k$ for all $\lambda \notin ER(p_k)$, in particular for $\nu$.  
%By our assumptions, $\alpha \notin \Phi^+_{max}$, so $Y^\alpha_k \equiv 0$.
Let $x_k = e^{\xi_k}$ for $\xi_k \rightarrow 0$ in $\lieg_{-\alpha}$.  Then, for some $m \in \NN$,
 $$ \Ad(x_k^{-1}) Y_k = Y_k^{\prime}=Y_k +\sum_{j=1}^{m} \frac{(-1)^j}{j!}(\ad \ \xi_k)^j(Y_k)$$
 %% Because $Y_k^{\alpha}=0$, we have $Y_k^{\prime} \in \lien^+$ for every $k$,
%Observe that for each $\lambda \in \Phi_{max}^+$, $\alpha_{j_0}$ occurs with a positive coefficient in the decomposition of
%%  $\lambda - l \alpha$,
%%   for every $l \in \NN$.   It follows that if $\lambda - l \alpha$ is a root, it is a positive root, and is in $(\Lambda^+)^c$ since $\alpha_{j_0}$
%%    occurs in its decomposition into simple roots.  
By our hypotheses, $Y_k^{\prime} \in \lien_\Lambda^+$,  hence $n_k' = e^{Y_k^{\prime}} \in P$ and $a_k n_k' \in H_{\Lambda}$.  
By Proposition \ref{prop.minus.alpha.conjugation}, we can choose $\xi_k \rightarrow 0$ in $\lieg_{- \alpha}$ such that the sequence $(Y_k^{' \nu})$ is unbounded.   

We have the relation 
$$a_kn_ke^{\xi_k}=e^{\Ad(a_k) \xi_k}a_k n_k'.$$ 
We wish to show that $\Ad(a_k) \xi_k \rightarrow 0$.
The action
     of $\Ad(a_k)$ on $\lieg_{- \alpha}$ is scalar multiplication by $\lambda_k=e^{- \alpha(Z_k)}$, where  $Z_k = \ln a_k$, so it is enough to show that
      $\lambda_k \leq C $, for some constant $C \in \RR$.  If this were not the case, then, up to taking a subsequence, there would be 
       $\zeta_k \rightarrow 0$ in $\lieg_{- \alpha}$ with $\Ad(a_k)\zeta_k \to \zeta_{\infty} \neq 0$.  For the product        
$$ p_ke^{\zeta_k}=e^{\Ad(a_k)\zeta_k} a_k e^{- \zeta_k} n_k e^{\zeta_k}$$
we know from above that $a_k e^{-\zeta_k} n_k e^{\zeta_k} \in P$.  Thus $p_k.[\zeta_k]=[\Ad(a_k)\zeta_k] \rightarrow [\zeta_\infty]$, while $L([\zeta_k]) \rightarrow 0$, which contradicts the fact that $(p_k)$ acts equicontinuously with respect to segments.

%The assumption on $\alpha(Z_k)$ implies $\Ad(a_k)$ is bounded in $\mbox{End}(\lieg_{-\alpha})$.  It follows that $\lim_{k \to \infty} \sup_{s \in [0,1]} (\Ad a_k)(s {\xi_k})=0$, hence 
Now let $\eta_k = \Ad(a_k) \xi_k$, which tends to 0.  It is easy to verify that 
$$ e^{-s \eta_k} p_k e^{s \xi_k} \in P \qquad \forall \ s \in \RR$$ 
Thus $q_k=a_k n_k'$ is a transverse perturbation of $(p_k)$ according to Definition \ref{defi.transverse.perturbation}, and, because $(Y_k'^\nu)$ is unbounded, it has $\nu \in ER(q_k)$, as desired. 
\end{preuve}

\begin{proposition}(Vertical sliding)
\label{prop.vertical.sliding}
Let $\nu \in (\Lambda^+)^c$ and $\alpha \in \Lambda^+$.  Let $p_k = a_k n_k \in H_\Lambda$ with $\alpha(Z_k) \geq M > - \infty$ ($\alpha(Z_k) \leq M < \infty$).  If $\nu + \alpha \in ER(p_k)$ (or $\nu - \alpha \in ER(p_k)$, resp.), then left and right vertical perturbation of $(p_k)$ yields $q_k = a_k n_k' \in H_\Lambda$ such that $\nu \in ER(q_k)$.
\end{proposition}

\begin{preuve}
%The proof is similar to the previous one.  
We can assume after vertical perturbation that $Y_k^\nu \equiv 0$.  We apply proposition \ref{prop.minus.alpha.conjugation} to obtain $\xi_k \rightarrow 0$ in $\lieg_{- \alpha}$ such that $(Y_k^{' \nu})$ is unbounded, where 
$$ Y_k^{\prime}= \Ad(x_k^{-1}) Y_k = Y_k +\sum_{j=1}^{m} \frac{(-1)^j}{j!}(\ad \ \xi_k)^j(Y_k)$$
for some $m \in \NN$, with $x_k = e^{\xi_k}$.  In this case, $Y_k \in \lien_\Lambda^+$ and $\alpha \in \Lambda^+$ together imply that $(\ad \xi_k)^j(Y_k) \in \lien_\Lambda^+$ for all $j \in \NN$.  Thus $Y_k' \in \lien_\Lambda^+$.

Let $n_k' = e^{Y_k'}$.  The lower bound on $\alpha(Z_k)$ implies $(\Ad a_k)(\xi_k)  \rightarrow 0$, so
$$e^{- \Ad(a_k)\xi_k} a_k n_ke^{\xi_k}=a_k n_k'$$ 
is obtained by left and right vertical perturbation from $(p_k)$. 

The proof for $\alpha(Z_k) \leq M < \infty$ and $Y_k^{\nu - \alpha}$ unbounded is similar.
\end{preuve}

\subsection{Algebraic proposition to reduce rank}

Using the tools developed so far in this section, we will now state the algebraic proposition that drives our induction step.  The next section contains the geometric interpretation of this result, and explains how to prove Theorem \ref{thm.induction} by induction on $\rk_\RR G$.  

\begin{proposition}
\label{prop.lower.rank}
Let $(p_k) = (a_k n_k)$ be a sequence of $H_\Lambda$ with $(n_k)$ unbounded. Assume that $(p_k)$, together with all its admissible perturbations, acts
 equicontinuously with respect to segments. Then an admissible perturbation of $(p_k)$ yields $(q_k)$ such that $ER(q_k)$ contains a root in $(\Lambda^+)^c \backslash \Phi^+_{max}$.
\end{proposition}

The proof of this proposition is given in Sections \ref{sec.proof.alg.prop} and \ref{sec.proof.alg.prop.g2} below.
   
%   
%   Assume first that we performe a Weyl reflection $r_{\alpha}$ on $p_k$, with $\alpha \in \Lambda$.  We already noticed that 
%   $q_k=r_{\alpha}(p_k)$
%    is still in $H_{\Lambda}$, and since $r_{\alpha}$ corresponds to an inner automorphism, there exists $p_{\alpha} \in G$ such
%     that  
%    $q_k=p_{\alpha}p_k p_{\alpha}^{-1}$.  Now, $Ad(p_{\alpha})$ leaves $\liep$ invariant, so that $p_{\alpha} \in Nor(P^{\,o})$, hence 
%     $p_{\alpha} \in P$ by our assumptions.  Now $\hx_k.p_{\alpha^{-1}}$ is bounded in $\hm$ and still project on $x_k$.  The sequence
%      $f_k(\hx_k.p_{\alpha}^{-1}).q_k^{-1}$, which is nothing else than $f_k(\hx_k).p_k^{-1}p_{\alpha}^{-1}$, is still bounded.
%       It follows that $q_k$ is an holonomy sequence at $x$, which is moreover a sequence in $H_{\Lambda}$, hence a reduced holonomy sequence. 
%       %We can then write $f_k(\hx_k.p_{\alpha})=$
%       
  
\section{Proof of Theorem \ref{thm.induction} by induction on rank}
\label{sec.induction}
The first half of this section gives the proof of Theorem \ref{thm.induction} from Proposition \ref{prop.lower.rank}.  The second half gives the proof of Proposition \ref{prop.lower.rank}.

\subsection{Invariant parabolic subvarieties}
Let $X=G/P$ with $G$ semisimple of real-rank $r$ and $P$ a parabolic subgroup with
 a Lie algebra $\liep= \liep_{\Lambda}$, $\Lambda \subsetneq \Phi$.
Let $V \subset X$ be a parabolic subvariety through the base point $o$. (These will be defined precisely below.)
%  $V=G^{\prime}/Q$ with $G^{\prime} \subset G$.  
If $(p_k)$ acts equicontinuously with respect to segments on $X$ and preserves $V$, then clearly it is equicontinuous with respect to segments on $V$.  The strategy for our induction argument is to find $(p_k)$-invariant $V \subset X$ of rank less than $r$.

Recall the notation introduced in Section \ref{sec.notation}, and denote by $B$ the Killing form on $\lieg$.
Given a subset $\Psi \subset \Phi$, 
%generates a semisimple root system $\Psi^+ \cup - \Psi^+$**.  
 let $\liea_0$ and $\liem_0$ be the ideals of 
 $\liea$ and $\liem$, respectively, commuting with $\oplus_{\alpha \in {\Psi}^+} \lieg_{\alpha}$.  Let $\liea_{\Psi}=\liea_0^{\perp}$ and $\liem_{\Psi}=\liem_0^{\perp}$, where the orthogonal is taken with respect to the scalar product $\langle X,Y \rangle =-B(X,\Theta Y )$.
 We obtain a subalgebra of $\lieg$ 
$$\lieg_{\Psi} = \sum_{\alpha \in \Psi^-} \lieg_{\alpha} \oplus \liea_\Psi \oplus \liem_\Psi \oplus \sum_{\alpha \in \Psi^+} \lieg_\alpha.$$
  It is easy to check that $\lieg_{\Psi}$ is $\Theta$-invariant, hence reductive,
    and has trivial center.  It follows that $\lieg_{\Psi}$ is semisimple.

%where $\liem_\Psi$ is a compact subalgebra centralizing $\liea_\Psi$ and contained in $\liem$.**

The corresponding connected subgroup $G_\Psi < G$ is closed. Indeed, $\ad (\lieg_{\Psi})$ is a semisimple subalgebra of 
$\End(\lieg)$, hence is an algebraic subalgebra (see \cite[Th 3.2]{hochschild}). For $G_{\Psi}^{\prime}$  the corresponding 
Zariski closed subgroup of $\GL(\lieg)$, the group $\Ad^{-1}(G_{\Psi}^{\prime})$ is closed in $G$, and so is its 
identity component $G_\Psi$.
%with parabolic subalgebra $\liep_{\Psi}$ as in subsection \ref{sec.notation.reduced}. $\lieh_{\Psi}$  

A minimal parabolic of $G_\Psi$ is contained in $P_{min}$.
The stabilizer of $o$ in $G_{\Psi}$ contains $P_{min} \cap G_{\Psi}$ and is algebraic, hence is a parabolic subgroup of $G_{\Psi}$, denoted $Q_{\Psi}$.  The orbit $G_{\Psi}.o$ is a \emph{parabolic subvariety} $V_{\Psi} \cong G_{\Psi}/Q_{\Psi}$, nontrivial provided $\Psi \not \subset \Lambda$, and of rank less than $r$.

 \begin{proposition}
 \label{prop.invariance}
Let $p_k = a_k n_k \in H_\Lambda$ and let $((Z_k^i),(Y^\alpha_k))$ be the exponential coordinates of $p_k$.
Then for any $\Psi \subset \Phi$, the variety $V_{\Psi} \subset X$ is invariant by $(p_k)$.  If $Z_k^i = 0$ for all $\alpha_i \in \Psi$, then $a_k$ acts trivially on $V_\Psi$; if $Y^\alpha_k = 0$ for all $\alpha \in \Psi^+ \cap (\Lambda^+)^c$, then $n_k$ is trivial on $V_{\Psi}$.
\end{proposition}

\begin{preuve}
Let $\xi  \in \sum \limits_{\alpha \in {\Psi^+}} \lieg_{-\alpha}$ and $x = e^\xi$.

Given $(Z_k)$ as in the hypotheses above, $\alpha(Z_k) \equiv 0$, for all $\alpha \in \Psi^+$.  Thus $\ad(\xi) Z_k =0$ and $\Ad(x) Z_k = Z_k$ for all $k$.  Thus $a_k x.o = x a_k.o = x.o$, and $a_k$ acts trivially on $V_\Psi$.

Now let $Y \in \lien_\Lambda^+$ with $Y^\alpha = 0$ for all $\alpha \in \Psi^+$.  
%Assume that $Y \in \Sigma_{\alpha \in (\Psi^+)^c} \lieg_{\alpha}$. 
  Write
 $$ \Ad(x) Y = Y' = Y + \sum_{k=1}^{m} \frac{(-1)^k}{k!}(\ad \xi)^{k}(Y) $$
 Note that $Y^{' \lambda} = 0$ unless $\lambda=\mu+\nu$, with $\mu$ a sum with negative integral coefficients of elements of $\Psi$ and $\nu$ in 
  $(\Psi^+)^c$; in particular, $\mu+\nu$ has positive coefficient on some simple root of $\Phi \backslash \Psi$.  In this case, $\lambda$ is a positive root, so $Y' \in \lien^+,$ and $e^{Y'} \in P$.  Thus $e^Yx.o= xe^{Y'}.o = x.o$, and $e^Y$ is trivial on $V_{\Psi}$.
   
The above calculation with $Y \in \sum \limits_{\alpha \in \Psi^+} \lieg_{\alpha}$ shows that $V_{\Psi}$ is invariant by $e^{Y}$; it is easy to see that $A$ leaves $V_\Psi$ invariant.  For invariance under a general sequence $p_k=a_kn_k$ in $H_\Lambda$, we can use the following basic lemma, the proof of which we leave to the reader:
   \begin{lemme}
   \label{lem.produit}
    Let $N$ be a simply connected nilpotent Lie group with Lie algebra $\lien$.  Let $\lien_0$ be an ideal of $\lien$, and let $Y, Y_0$ be elements of
     $\lien$ and $\lien_0$.  Then there exists $Y_0^{\prime} \in \lien_0$ such that
     $$ e^{Y+Y_0}=e^Ye^{Y_0^{\prime}}.$$
   \end{lemme}

   This lemma allows to write $n_k=e^{W_k}e^{U_k}$ with $W_k \in  \sum \limits_{\alpha \in \Psi^+} \lieg_{\alpha}$ and
   $U_k \in \sum \limits_{(\Psi^+)^c} \lieg_{\alpha}$.  We can then conclude because   each factor $a_k$, $e^{W_k}$ and $e^{U_k}$ preserves $V_{\Psi}$.
   \end{preuve}

The unipotent radical of $Q_\Psi$ is $N_{\Psi,\Lambda}^+ < N_\Lambda^+$ with Lie algebra
$$\lien_{\Psi,\Lambda}^+ = \oplus_{\alpha \in \Psi^+ \backslash \Lambda^+} \lieg_{\alpha}$$
The analogue of $H_\Lambda$ in $G_\Psi$ is
$H_{\Psi,\Lambda} = A_\Psi \ltimes N_{\Psi,\Lambda}^+$.  Note that 
$$N_\Lambda^+ = N^+_{\Psi,\Lambda} \cdot (N^+_\Psi \cap N_\Lambda^+),$$ 
and that the second factor is normal in $H_\Lambda$.  We will also need below the decomposition $A = A_\Psi \cdot A_{\Phi \backslash \Psi}$. 
%, where $A_{\Psi \cap \Lambda}$ is the connected subgroup of $A$ with Lie algebra spanned by the coroots dual to the roots in $\Psi \cap \Lambda$.

\subsection{The induction step}

Suppose that Theorem \ref{thm.induction} holds for all parabolic models $G/P$ of real-rank at most $r-1$.  
We will prove using Proposition \ref{prop.lower.rank} that it holds for all models of real-rank $r$.  
Let $X = G/P_\Lambda$ of rank $r$ be given, and let $(p_k)$ be a sequence of $H_\Lambda$ which, together with all its admissible perturbations, acts
 equicontinuously with respect to segments.  
 The aim is to show that $(n_k)$ is bounded.  If not, then Proposition \ref{prop.lower.rank} gives, after an admissible perturbation, $(q_k)$ with $ER(q_k)$ containing a root $\lambda \in (\Lambda^+)^c \backslash \Phi^+_{max}$.  

There is a proper subset $\Psi$ of $\Phi$ such that $\lambda \in \Psi^+$.  It cannot be that $\Psi$ is contained in $\Lambda$ because $\lambda \in (\Lambda^+)^c$. Now $q_k \in H_\Lambda$ preserves $V_\Psi$ by Proposition \ref{prop.invariance}; denote the restriction by $(q'_k)$, which is a sequence of $Q_\Psi$, and let $a_k'n_k'$ be the decomposition into components on $A_{\Psi}$ and $N^+_{\Psi,\Lambda}$, respectively. Because $\lambda \in ER(q_k)$, it follows that $(n_k')$ is unbounded.

As $\rk_\RR G_\Psi \leq r-1$, the induction hypothesis yields a contradiction, \emph{provided that all admissible perturbations of $(q_k')$ in $G_\Psi$ act equicontinuously with respect to segments on $V_\Psi$}.  Admissible perturbation in $G_\Psi$ means more precisely that vertical and transverse perturbations are as in Section \ref{sec.vertical.transverse} with $\lieg_{\Psi}$ in place of $\lieg$, and $Q_{\Psi}$ in place of $P$, and Weyl reflections are done with respect to roots $\alpha$ in $(\Psi \cap \Lambda)^+$. The following lemma ensures that $(q_k')$ satisfies the hypotheses of Theorem \ref{thm.induction} and allows us to apply our induction hypothesis: 

 \begin{lemme}
 \label{lem.induction}
 Let $X=G/P_\Lambda$ be a parabolic variety, and $(q_k)$ be a sequence of $H_{\Lambda}$.  Assume that $(q_k)$ preserves a parabolic subvariety $V_{\Psi}$ on which it restricts to $(q_k')$.  If every admissible perturbation of $(q_k)$ acts equicontinuously with respect to segments in $X$, 
  then every admissible perturbation of $(q_k')$ in $G_{\Psi}$ acts equicontinuously with respect to segments in $V_{\Psi}$.
  \end{lemme}

\begin{preuve}
    We will prove that any admissible perturbation of the sequence $(q_k')$ in $G_{\Psi}$ can be obtained by an admissible perturbation of $(q_k)$, restricted to $V_{\Psi}$.  Assume that $(p_k')$ is obtained from $(q_k')$ by an admissible perturbation in $G_{\Psi}$. 
  We seek an admissible perturbation
   $(p_k)$ of $(q_k)$, such that $p_k$ preserves $V_{\Psi}$, and the restriction of $p_k$ to $V_{\Psi}$ is precisely 
   $p_k^{\prime}$.  Existence of such $(p_k)$ can be checked for each of the three kinds of admissible perturbations in $G_\Psi$:  
%The notations used below are those of Section    \ref{sec.invariant.varieties}
   
%   \begin{itemize}
{\bf (1) vertical perturbation}: There are bounded sequences $(l_k)$ and $(m_k)$ in $Q_{\Psi}$ such that 
    $p_k^{\prime}=l_kq_k'm_k$ on $V_\Psi$.  Because  $Q_{\Psi} < P$, the desired vertical perturbation of $(q_k)$ in $G$ is simply $(p_k) = (l_kq_km_k)$.

\smallskip

{\bf (2) transverse perturbation}: In this case, write
     $p_k^{\prime}=e^{- \eta_k}q_k'e^{\xi_k}$ where  $(\eta_k)$ and $(\xi_k)$ are two sequences of $\lieg_{\Psi} \setminus {\mathfrak q}_{\Psi}$ tending to $0$.  As these are also sequences of $\lieg \setminus \liep$, we can set $p_k=e^{- \eta_k}q_ke^{\xi_k}$; we will show that this is a transverse perturbation in $G$.  

Let $x \in V_{\Psi}$. Observe that because 
       $\xi_k, \eta_k \in \lieg_{\Psi}$,
       $$e^{- s\eta_k}q_ke^{s\xi_k}.x=e^{- s\eta_k}q_k'e^{s\xi_k}.x \qquad \forall \ s \in \RR;$$
thus $e^{- s\eta_k}q_ke^{s\xi_k}$ preserves $V_{\Psi}$ and acts on it by $e^{-s \eta_k}q_k'e^{s\xi_k}$.
        Taking 
         $x=o$ gives $e^{- s\eta_k}q_ke^{s\xi_k}.o=e^{- s\eta_k}q_k'e^{s\xi_k}.o=o$, because the latter is
         in $Q_\Psi$ for all $s$. This proves  $e^{- s\eta_k}q_ke^{s\xi_k} \in P$ for all $s \in \RR$, and $p_k$ is a transverse perturbation of
        $q_k$.
    
\smallskip

{\bf (3) Weyl reflection}: Let $r_{\alpha} \in \Aut(G_\Psi)$ realize the Weyl reflection $\rho_\alpha$, for $\alpha \in (\Psi \cap \Lambda)^+$. 
    Decompose, using Lemma \ref{lem.produit}, 
       $$q_k=a_kn_k = a_k^{\prime \prime}a_k^{\prime}n_k^{\prime}n_k^{\prime \prime},$$ 
where $a_k^{\prime} \in A_{\Psi}$, $n_k^{\prime} \in N_{\Psi,\Lambda}^+$, $a_k^{\prime \prime} \in A_{\Phi \backslash \Psi}$, and $n_k^{\prime \prime} \in (N_\Psi^+ \cap N^+_\Lambda)$. 
By Proposition \ref{prop.invariance}, both $a_k^{\prime \prime}$ and $n_k^{\prime \prime}$ are in the kernel of the restriction to $V_\Psi$, so we can write $q_k'=a_k^{\prime}n_k^{\prime}$.
        
        Now let
     ${\tilde r}_{\alpha}$ be an automorphism of $G$ effecting $\rho_{\alpha}$ 
     on $\liea^*$.  Because $\alpha \in (\Psi \cap \Lambda)^+$,
      the derivative of $\tilde{r}_{\alpha}$ preserves the Lie algebras $\liea_{\Psi}$, $\liea_{\Phi \backslash \Psi}$, $\lien_{\Psi,\Lambda}^+$ and $(\lien_{\Psi}^+ \cap \lien_{\Lambda}^+)$, so $\tilde{r}_\alpha$ preserves the corresponding connected subgroups in $G$.  Thus $\tilde{r}_{\alpha}(q_k')=r_{\alpha}(q_k')$, and 
        $$\tilde{r}_{\alpha}(q_k)=\tilde{r}_{\alpha}(a_k^{\prime \prime})r_{\alpha}(q_k')\tilde{r}_{\alpha}(n_k^{\prime \prime})$$
      %  Since $\tilde{r}_{\alpha}(a_k^{\prime \prime})$ and $\tilde{r}_{\alpha}(n_k^{\prime \prime})$ belong respectively to $A_0$ and $ (N_{\Psi,\Lambda}^+)^c$, 
preserves $V_{\Psi}$
%and act trivially on it (see Proposition \ref{}).  Thus $\tilde{r}_{\alpha}(p_k)$ preserves         $V_{\Psi}$ and 
and restricts on it to $r_{\alpha}(q_k')$, as desired.
        %        \end{itemize}
\end{preuve}

%% Now, for the proof of Theorem \ref{} in higher rank, we can state the 

%% {\bf Induction Hypothesis}:  \emph{Let $G'$ be any semisimple Lie group of rank less than $r$, and $P_\Lambda'$ a parabolic subgroup.  Let $(p_k) = (a_k n_k)$ be a sequence of $H_\Lambda$.  Given $V_\Psi \subset X$ a parabolic subvariety of rank less than $r$, denote $(a_k' n_k')$ the restriction of $(p_k)$ to $V_\Psi$, with $a'_k \in A_{\Psi \cap \Lambda}$ and $n_k' \in N^+_{\Psi,\Lambda}$.  Then $(n_k')$ is bounded.}

The proof by induction of Theorem \ref{thm.induction} is now complete, once we prove Proposition \ref{prop.lower.rank}.

\subsection{Proof of Proposition \ref{prop.lower.rank} (assuming the root system of $\lieg$ is not of type $G_2$)}
\label{sec.proof.alg.prop}

%% \begin{proposition}
%%  \label{prop.technical}
%%   Let $G$ be simple of real rank $r \geq 2$, and let $P=P_{\Lambda}$. 
%%     Let $p_k=a_kn_k \in H_{\Lambda}$, with $(n_k)$ unbounded. 
%%    Then finitely-many admissible operations convert $(p_k)$ to $(q_k)$ in $H_{\Lambda}$ leaving invariant a parabolic subvariety $V_{\Psi}=G_{\Psi}/Q_{\Psi}$ with $\rk_\RR(G_\Phi) < \rk_\RR(G)$, and on which $q_k$ acts by $a_k^{'}n_k^{'} \in H_{\Psi, \Lambda}$ with $n_k'$ unbounded in $N_{\Psi, \Lambda}^+$.
%%  \end{proposition}

Proposition \ref{prop.lower.rank} is vacuously true if the set $\Phi^+_{max}$ is empty.  Thus, we assume from now on 
that \emph{$G$ is a simple Lie group}.

Let $(p_k) = (a_k n_k)$ be a sequence of $H_\Lambda$ with $(n_k)$ unbounded.  
%Assume that $(p_k)$ acts equicontinuously with respect to segments on $X = G/P$.
%% If $ER(p_k)$ contains a root $\lambda \not \in \Phi^+_{max}$, then $\lambda \in \Psi^+$ for some $\Psi \subsetneq \Phi$.  Then by Proposition \ref{prop.invariance}, $(p_k)$ restricts on $V_{\Psi}$ to a sequence $a_k^{'}n_k^{'}$, with $n_k'$ unbounded in $N_{\Psi,\Lambda}^+$, so conclusion (2) holds. 
That means $ER(p_k) \subseteq (\Lambda^+)^c$ is nonempty.  If it contains a root not in $\Phi^+_{max}$, then there is nothing to show, so we suppose that $ER(p_k) \subseteq \Phi^+_{max}$.  Define the \emph{degree} of $\alpha \in \Phi^+$ to be the sum of the coefficients in the unique expression of $\alpha$ as a positive integral linear combination of roots in $\Phi$.

 Let $Y_k = \ln n_k$.  By Proposition \ref{prop.trivial.or.unbounded}, we may assume $Y_k^\lambda \equiv 0$ for $\lambda \notin ER(p_k)$.  To prove that an admissible perturbation of $(p_k)$ results in $(q_k)$ with $ER(q_k)$ not contained in $\Phi^+_{max}$, we will show that 
%(\emph{a}) we can directly find such a perturbation; or (\emph{b}) 
for any $\lambda \in ER(p_k)$ of minimal degree, there is a sequence of admissible operations resulting in $\lambda' \in ER(q_k)$ with the degree of $\lambda'$ strictly lower than the degree of $\lambda$.

Let $\lambda \in ER(p_k) \subseteq \Phi^+_{max}$ of minimal degree.  There is some $\alpha \in \Phi$ with $\langle \alpha, \lambda \rangle > 0$; otherwise, $\lambda$ would be in the negative of the Weyl chamber spanned by $\Phi$, contradicting that it is a positive root.
  For such $\alpha$,
$$ A_{\alpha \lambda} = \frac{2 \langle \alpha, \lambda \rangle}{\langle \alpha, \alpha \rangle} > 0$$

\smallskip

{\bf Case $\alpha \in \Lambda$.}  In this case, the Weyl reflection $\rho_\alpha \in W_\Lambda$ yields
$$ \rho_\alpha (\lambda) = \lambda' = \lambda - A_{\alpha \lambda} \alpha \in (\Lambda^+)^c$$
of smaller degree.  The admissible operation $r_\alpha$ yields $q_k \in H_{\Lambda}$ with $\lambda' \in ER(q_k)$.

%% Let $\beta \in \Phi \backslash \Lambda$ be given by Proposition \ref{prop.weyl.transitivity} applied to $\alpha$.  By lemma \ref{lem.invariance}, $V_{\{ \beta \}}$ is $\mbox{Hol}_\Lambda$-invariant.  By the induction hypothesis, theorem \ref{thm.equicontinuity} holds for the restriction of $(p_k)$ to $V_{ \{ \beta \}}$.  In particular, $\beta(Z_k)$ is bounded away from $- \infty$.
%% %also, using proposition \ref{prop.trivial.or.unbounded}, we may assume $Y_k^\beta \equiv 1$. 

%% Now proposition \ref{prop.vertical.sliding} applies to $(p_k)$, with $\nu  = \lambda_1 = \lambda_0 - \alpha$, to give a perturbation $(q_k = a_k n_k')$ of $(p_k)$ such that $Y_k^{' \lambda_1}$ is unbounded.  Applying proposition \ref{prop.vertical.sliding} $l-1$ more times results in a holonomy sequence $(q_k = a_k n_k')$ with $Y_k'$ unbounded along $\lieg_{\lambda_l}$, as promised.

\smallskip

{\bf Case $\alpha \in \Phi \backslash \Lambda$.}  Note that 
$\nu = \lambda - \alpha \in \Phi^+$, because $\lambda - A_{\alpha \lambda} \alpha \in \Phi^+$, and strings are unbroken.  

If $P = P_\Lambda$ is not a maximal parabolic with $\Lambda = \Phi \backslash \{ \alpha \}$, then $(p_k)$, $\alpha$, and $\nu$ satisfy the hypotheses of Proposition \ref{prop.transverse.sliding}, which  
%in Remark \ref{rmk.admissible.roots}, so $- \alpha$ is admissible for $(p_k)$.
%Propositions \ref{prop.simple.transverse.sliding} and \ref{prop.transverse.sliding} give a sliding 
thus gives another holonomy sequence 
 $(q_k)$ with $\nu = \lambda - \alpha \in ER(q_k)$, which has lower degree than $\lambda$. 

% 
% By lemma \ref{lem.invariance}, $V_{\{ \alpha \}}$ is $H_\Lambda$-invariant.  By the induction hypothesis, theorem
% \ref{thm.equicontinuity} holds for the restriction of $(p_k)$ to $V_{ \{ \alpha \}}$. 
% Using proposition \ref{prop.trivial.or.unbounded}, we may assume $Y_k^\alpha \equiv 1$; 
% we also have $\alpha(Z_k)$ bounded away from $- \infty$.  
% Now apply proposition 
% It also belongs to $(\Lambda^+)^c$ unless

Now suppose $P$ is a maximal parabolic, with $\Lambda = \Phi \backslash \{ \alpha \}$.    Every root in $ER(p_k)$ has the form $\lambda_i = m_i \alpha + \mu_i$, where $m_i \geq 1$, and $\mu_i$ is in the positive integral span of $\Lambda$.  If none of the $\mu_i$ is a root, then again the hypotheses of Proposition \ref{prop.transverse.sliding} are satisfied,
%is admissible, 
so, as above, there is a holonomy sequence $(q_k)$ with $\lambda - \alpha \in ER(q_k)$.

Thus we may assume that $\mu_i$ is a root for some $i$.  

\begin{lemme}
Let $P_\Lambda < G$ be a maximal parabolic with $\Lambda = \Phi \backslash \{ \alpha \}$.  If $m \alpha + \mu \in \Phi^+_{max}$ for $m \geq 1$ and $\mu \in \Lambda^+$, then $\alpha$ is a valence-one vertex of the Dynkin graph of $\lieg$---that is, $A_{\alpha \beta} \neq 0$ for exactly one element $\beta \in \Lambda$.
\end{lemme}

\begin{preuve}
The positive root $\mu$ can be written $\alpha_1 + \cdots + \alpha_\ell$, with $\alpha_i \in \Phi$ and  $\alpha_1 + \cdots + \alpha_i \in \Phi^+$ for all $i$.  This statement is established by induction on the degree of $\mu$ together with the fact, seen above, that there exists at least one $\alpha \in \Phi$ with $\langle \alpha, \mu \rangle > 0$.  On the other hand, the expression in terms of simple roots is unique, and $\mu \in \Lambda^+$, so $\alpha_i \in \Lambda$ and $\alpha_1 + \cdots + \alpha_i \in \Lambda^+$ for all $i$.  Denote this latter root $\mu_i$.

Next we show by induction that the collection of roots $\{ \alpha_1 , \ldots, \alpha_i \}$ corresponds to a connected set of vertices in the Dynkin diagram of $\lieg$ for all $i$.  Assume connectedness for a given $i \geq 1$.  If $A_{\mu_i \alpha_{i+1}} = 0$, then the $\alpha_{i+1}$-string of $\mu_i$ is symmetric, of the form $\{ \mu_i - q \alpha_{i+1}, \ldots, \mu_i, \ldots, \mu_i + q \alpha_{i+1} \}$.  Note that $q \geq 1$ because $\mu_{i+1} \in \Lambda^+$.  Then $\mu_i - \alpha_{i+1} = \lambda$ is also a positive root, with $\mu_i = \lambda + \alpha_{i+1}$.  In other words, $\alpha_{i+1}$ already belongs to the collection of roots appearing in $\mu_i$, which is connected by the induction hypothesis.
%The root $\mu$ belongs to some basis of simple roots, and the Weyl group $W$ acts transitively on such sets (see \cite[Th 2.6.3] 
%%p 119]
%{knapp.lie.groups}), which means there is $\rho \in W$ sending some $\alpha_i \in \Phi$ to $\mu$.  This $\rho$ is moreover a product $\rho_{i_\ell} \cdots \rho_{i_1}$ of Weyl reflections.   Let $\mu_0 = \alpha_i$ and $\mu_j$ be the result after performing $j$ reflections.  Then one can see that at each step, $\mu_j$ is a positive root, comprised of simple roots that form a connected subset of the Dynkin graph.  If $\rho_{i_j}$ is the reflection at the $j$th step, then $\alpha_{i_j}$ is connected to exactly one of the simple roots appearing in $\mu_{j-1}$ because the Dynkin diagram is a tree, and it adds a positive multiple of $\alpha_{i_j}$ to make $\mu_j$.  
%

We conclude that the elements of $\Phi$ appearing in the decomposition of $\mu$ correspond to a connected subset of the Dynkin graph.  These are precisely the elements of $\Lambda = \Phi \backslash \{\alpha \}$.  As the Dynkin graph is a connected tree, the conclusion follows.  
\end{preuve}

Let $\beta \in \Lambda$ with $A_{\alpha \beta} \neq 0$.  Write $\lambda_i = \lambda' = m' \alpha + \mu'$ where $\mu' \in \Lambda^+_{max}$, and let $c' \in \BZ^+$ be the coefficient of $\beta$ in $\mu'$.  
%Now $A_{\alpha \mu'} = c_1' A_{\alpha \beta}$ and $A_{\mu' \alpha} = c_1' A_{\beta \alpha}$.  
The product
$$ A_{\alpha \mu'} A_{\mu' \alpha} = \frac{(c')^2 A_{\alpha \beta} A_{\beta \alpha} \langle \beta, \beta \rangle}{\langle \mu', \mu' \rangle} \in \{ 1,2,3 \}.$$

(Although our root system is not necessarily reduced, the value 4 could only occur for $\mu' = 2 \alpha$ or $\alpha = 2 \mu'$, neither of which is the case.)
First suppose the Dynkin diagram has no double or triple edges, so the root system of $\lieg$ is $A_r, D_r, E_6, E_7,$ or $E_8$.  Then all roots of $\Phi^+$ have the same length and $A_{\alpha \beta} A_{\beta \alpha} = 1$.  In this case, $A_{\alpha \mu'} A_{\mu' \alpha} = (c')^2$, so $c'=1$ and $A_{\alpha \mu'} = -1 =  A_{\mu' \alpha}$.  The $\alpha$-string of $\mu'$ comprises $\mu'$ and $\mu' + \alpha$.  Hence $m'=1$ and $\lambda' = \mu' + \alpha$.  The $\mu'$-string of $\alpha$ comprises $\alpha, \lambda'$.  Now $\rho_{\mu'}(\lambda') = \alpha$, so the Weyl reflection $r_{\mu'}(p_k)$ is an admissible perturbation resulting in $(q_k)$ with $\alpha \in ER(q_k)$. 

Under the assumption that $\lieg$ is not of type $G_2$, there are no triple bonds in the Dynkin diagram of $\lieg$, so it remains to consider the root systems with double bonds: $B_r, BC_r, C_r$, and $F_4$.  Let $\lambda$ with $A_{\alpha \lambda} > 0$ as above be of minimal degree in $ER(p_k).$  Write $\lambda = m \alpha + \mu$, where $\mu$---not necessarily a root---is a positive integral combination of elements of $\Lambda$, and let $c \in \BZ^+$ be the coefficient of $\beta$ in $\mu$. 
Because $\lambda - A_{\alpha \lambda} \alpha$ is a positive root, 
\begin{equation}
\label{eqn.A.alpha.lambda}
0 < A_{\alpha \lambda} = 2m + c A_{\alpha \beta} \leq m
\end{equation}

Write $\Phi = \{ \gamma_1, \ldots, \gamma_r \}$, numbered from left to right in the Dynkin diagram, where we follow the ordering of \cite{knapp.lie.groups}.  We have $\alpha = \gamma_1$ or $\gamma_r$.

{\bf Type $B_r$ or $BC_r$}. 
%Next consider $A_{\alpha \beta} = -2$ and $A_{\beta \alpha} = -1$.  
For $B_r$, the set $\Phi^+_{max}$ comprises, for $i = 2, \ldots, r$,
$$ \lambda_1 = \gamma_1 +  \cdots + \gamma_r,  \ \lambda_i = \lambda_1 + \gamma_i + \cdots + \gamma_r,$$

If $\alpha$ is the short root, $\gamma_r$, then $A_{\alpha \beta} =-2$. The possibility $m=1$ is incompatible with (\ref{eqn.A.alpha.lambda}).  If $m=2$, then the same inequality implies $c = 1$, so $\lambda = \lambda_r$.  
If $r > 2$, then $\rho_{\gamma_1} (\lambda)$ has lower degree,
% = 2 \alpha + \beta_1 + \cdots + \beta_{r-2}$, 
so a Weyl reflection $r_{\gamma_1}$ is an admissible perturbation with the desired effect. 
Otherwise, $r=2$ and $\lambda = \beta + 2 \alpha$.  In this case, as $\lambda$ is an element of $ER(p_k)$ of minimal degree, $ER(p_k) = \{ \lambda \}$.  
There is a rank-one subvariety $V_\lambda \subset X$ left invariant by $(p_k)$ and on which it restricts to $(a_k' n_k)$ with $(n_k)$ unbounded.  Proposition \ref{prop.nonequicontinuous.rk1} leads to a contradiction.

If $\alpha$ is the long root $\gamma_1$, then $m=1$ and $A_{\alpha \beta} = -1$, so (\ref{eqn.A.alpha.lambda}) implies $c = 1$.  Then $\lambda = \lambda_1$ or $r > 2$ .   In the first case, $\mu = \gamma_2 + \cdots + \gamma_r \in \Lambda^+$ is a short root with $A_{\mu \alpha} = -2$.  Proposition \ref{prop.vertical.sliding} permits vertical sliding along $- \mu$, resulting in $(q_k)$ with $\alpha \in ER(q_k)$, or along $\mu$, resulting in $(q_k)$ with $\alpha + 2 \mu \in ER(q_k)$.  In the latter case, the Weyl reflection $\rho_{\mu}(\alpha + 2 \mu) = \alpha$, so $r_{\mu}$ leads to the desired conclusion. Otherwise, $\lambda = \lambda_i$ for $2 < i \leq r$; in this case, Weyl reflection in $\gamma_i \in \Lambda$ results in $(q_k)$ with a minimal element of $ER(q_k)$ of lower degree.

In $BC_r$, the set $\Phi^+_{max}$ comprises $\{ \lambda_i \ : \ 1 \leq i \leq r \}$ from above, together with $2 \lambda_1$.   If $\lambda = 2 \lambda_1$, then $ER(p_k) = \{ \lambda \}$; in this case, 
restricting to the rank-one subvariety $V_{\lambda}$ yields a contradiction to Proposition \ref{prop.nonequicontinuous.rk1}. 

\smallskip

{\bf Type $C_r$}. The set $\Phi^+_{max}$ comprises, for $i = 1, \ldots, r-1$,
$$\lambda_r =  \gamma_1 + \cdots + \gamma_r, \ \lambda_i = \lambda_r + \gamma_i + \cdots + \gamma_{r-1}$$

If $\alpha$ equals the long root, $\gamma_r$, then $A_{\alpha \beta} = -1$ and $m=1$.  The inequality (\ref{eqn.A.alpha.lambda}) gives $c = 1$ and $\lambda = \lambda_r$. If $r > 2$, then $A_{\gamma_1 \lambda} = 1$, and  $\rho_{\gamma_1}(\lambda)$ is a root of lower degree.  The remaining possibility is $r=2$ with $ER(p_k) = \{ \alpha + \beta, \alpha + 2 \beta \}$ or simply $\{ \alpha + \beta \}$.  In the first case, the Weyl reflection $r_\beta$ results in $(q_k)$ with $\alpha \in ER(q_k)$.  In the second case, we again apply Proposition \ref{prop.nonequicontinuous.rk1}. 

When $\alpha$ equals the short root $\gamma_1$, then we first consider $\lambda = \lambda_i$ for $i \neq 1$.  The Weyl reflection $\rho_{\gamma_i}(\lambda)$ has lower degree.  If $\lambda = \lambda_1$, then $ER(p_k) = \{ \lambda \}$, so Proposition \ref{prop.nonequicontinuous.rk1} completes the proof.

\smallskip

{\bf Type $F_4$.}  The roots in $\Phi^+_{max}$, in terms of the basis $\{ \gamma_i \}$, are 
\begin{eqnarray*}
(1,1,1,1), (1,1,2,1), (1,1,2,2), (1,2,2,1), (1,2,2,2), \\
(1,2,3,1), (1,2,3,2), (1,2,4,2), (1,3,4,2), (2,3,4,2)
\end{eqnarray*}

%In this case, $\alpha$ could be $\gamma_1$ or $\gamma_4$, which we number so that $\gamma_1$ is a long root and $\gamma_4$ a short one.  
Recall that $ER(p_k)$ contains $\lambda' = m' \alpha + \mu'$ with $\mu'$ a root in $\Lambda^+_{max}$.  The roots of $\Lambda^+_{max}$ correspond to those of $C_3$ when $\alpha$ equals the long root $\gamma_1$ and $B_3$ when $\alpha$ equals the short root $\gamma_4$.    In the first case the possibilities are
$$ \lambda' \in  \Lambda' = \{ (1,1,1,1), (1,1,2,1), (1,1,2,2) \}$$
%\gamma_1 + \gamma_2 + \gamma_3 + \gamma_4, \gamma_1 + \gamma_2 + 2 \gamma_3 + \gamma_4, \gamma_1 + \gamma_2 + 2 \gamma_3 +  2 \gamma_4 \}.$$

The maximum degree in $\Lambda'$ is 6.  As all other roots of $\Phi^+_{max}$ have degree at least 6, we may assume $\lambda \in ER(p_k)$ of minimal degree belongs to $\Lambda'$.  If $ER(p_k) = \{ \gamma_1 + \gamma_2 + 2 \gamma_3 +  2 \gamma_4  \}$, then we can invoke Proposition \ref{prop.nonequicontinuous.rk1}.  Otherwise, a Weyl reflection in $\gamma_4$ or $\gamma_3$ reduces the degree of $\lambda$ and yields a new holonomy sequence $(q_k)$ with an element of lower degree in $ER(q_k)$.

In the second case, $\Lambda'$ contains the roots listed above, together with 
$$ (1,2,2,1), (1,2,2,2) $$  
%\gamma_1 + 2 \gamma_2 + 2 \gamma_3 + \gamma_4, \gamma_1 + 2 \gamma_2 + 2 \gamma_3 + 2 \gamma_4$$
Now the maximal degree in $\Lambda'$ is 7, and all other roots of $\Phi^+_{max}$ have degree at least 7, so we may again assume $\lambda \in \Lambda'$.  A Weyl reflection in $\gamma_1$ or $\gamma_2$ will reduce the degree of any $\lambda \in \Lambda'$, giving the desired conclusion in this case.

\subsection{Proof of Proposition \ref{prop.lower.rank} for $G_2$}
\label{sec.proof.alg.prop.g2}

Assume $\lieg$ is of type $G_2$, and write $\Phi = \{ \alpha, \beta \}$ with $| \alpha | \leq |\beta|$.  Then
$$\Phi^+_{max} = \{ \alpha + \beta, 2 \alpha + \beta, 3 \alpha + \beta, 3 \alpha + 2 \beta \}$$

Assume first that $\Lambda = \{ \alpha \}$, so $A_{\alpha \beta} = -3$. Given $\lambda \in ER(p_k)$ of minimal degree, the goal is to find an admissible perturbation $(q_k)$ with $\beta \in ER(q_k)$.  As in the previous section (but with the roles of $\alpha$ and $\beta$ switched), we can assume that $A_{\beta \lambda} > 0$.  The two possibilities for $\lambda$ are thus $ 3\alpha + 2 \beta$ or $\alpha + \beta$.  In the first case, $\lambda$ is the only element of $ER(p_k)$, so we can conclude using Proposition \ref{prop.nonequicontinuous.rk1} as in the cases of $C_2$ and $B_2$.  In the second case, we apply Proposition \ref{prop.vertical.sliding}.  We can assume, after passing to a subsequence, that $\alpha(Z_k)$ is bounded either below or above.  If it is bounded below, then a vertical sliding on $(p_k)$ yields $(q_k)$ with $\beta \in ER(q_k)$, as desired.  If $\alpha(Z_k)$ is bounded above, then vertical slidings give $3 \alpha + \beta$ in $ER(q_k)$.  Then the Weyl reflection $r_\alpha$ on $(q_k)$ gives $(s_k)$ with $\beta \in ER(s_k)$.  

 %% so $\rho_{\alpha} (3 \alpha + \beta) = \beta$.  Thus, if any of $\alpha + \beta, 2 \alpha + \beta, 3 \alpha + \beta \in ER(p_k)$, then a combination of $s_{\alpha, \xi_k}$, $s_{- \alpha, \xi_k}$, according as $\alpha$ or $- \alpha$ is admissible, together with $r_\alpha$, will eventually yield $q_k$ with $\beta \in ER(q_k)$.

%% If $ER(q_k) = \{ 3 \alpha + 2 \beta \}$, then $- \beta$ is an admissible root, and a transverse sliding $q_k = s_{- \beta, \xi_k}(p_k)$ yields $3 \alpha + \beta \in ER(q_k)$.  Then $q_k' = r_\alpha(q_k)$ has $\beta \in ER(q_k')$.

Now consider $\Lambda = \{ \beta \}$, so $A_{\beta \alpha} = -1$.  The condition $A_{\alpha \lambda} > 0$ leaves the possibilities $2 \alpha + \beta$ or $3 \alpha + \beta$ for $\lambda$.  Unfortunately, the tools used above don't help in either of these cases.  The solution is to slide along $- \alpha$, although it does not satisfy the hypotheses of Proposition \ref{prop.transverse.sliding}. 

Let $S \cong Z(S) S_0$ be the reductive complement in a Levi decomposition of $P_\beta$, where $S_0$ is simple of rank one.  The group $S$ admits a $KAK$ decomposition, where $A = \exp(\liea)$ as defined above, and $K$ is a maximal compact subgroup of $S_0$. Write $N_\beta^+$ for the unipotent radical of $P_\beta$.  The decomposition of the corresponding Lie algebra $\lien_\beta^+$ into irreducible subspaces under $\Ad(S)$ is $E_1 \oplus E_2 \oplus E_3$, where $E_1 = \lieg_\alpha \oplus \lieg_{\alpha + \beta}$; $E_2 = \lieg_{2 \alpha + \beta}$; and $E_3 = \lieg_{3 \alpha + \beta} \oplus \lieg_{3 \alpha + 2 \beta}$.  This decomposition can be seen from the fact that $\lies$ is contained in the sum of root spaces $\lieg_{- \beta} \oplus \lieg_0 \oplus \lieg_{\beta}$. 

%\begin{definition}
%\label{defi.exceptional.sliding}
%Let $\Phi = \{ \alpha, \beta \}$ and $\Lambda = \{ \beta \}$ with $|\alpha| \leq |\beta|$.  Given a sequence $(p_k)$ in $H_\Lambda$, 

Recall that $p_k = a_k n_k$ with $Y^\nu_k \equiv 0$ if $\nu \notin ER(p_k)$.  Let $\xi_k \rightarrow 0$  in $\lieg_{- \alpha}$ and $x_k = e^{\xi_k}$, and set
$$ q_k = e^{- \Ad(a_k)\xi_k} p_k e^{\xi_k} = a_k x_k^{-1} n_k x_k.$$  
Just as in the proof of Proposition \ref{prop.transverse.sliding}, $\Ad(a_k) \xi_k \rightarrow 0$ and $(q_k)$ is a transverse perturbation of $(p_k)$; it is in particular a sequence in $P$, although it may not be in $H_\beta$.  More precisely, $x_k^{-1} n_k x_k \in N^+$, which can be deduced from the formula,
 $$ \Ad(x_k^{-1}) Y_k =Y_k +\sum_{j=1}^{m} \frac{(-1)^j}{j!}(\ad \ \xi_k)^j(Y_k)$$
with $Y_k = \ln n_k$.
 Using Lemma \ref{lem.produit}, write $q_k = a_k u_k n''_k$ with $a_k u_k \in S$ and $n_k'' \in N_\beta^+$.  Proposition \ref{prop.minus.alpha.conjugation} gives that $\lambda - \alpha \in ER(n_k'')$.  Performing this transverse sliding twice if necessary, depending on $\lambda$, we arrive at $\alpha + \beta \in ER(n_k'')$.
   
%contained in $P$ thanks to Proposition \ref{prop.simple.transverse.sliding}.  
Next, 
%write $q_k' = s_k n_k'$ according to the Levi decomposition of $P$, and 
let $l_k' a_k' l_k$ be the $KAK$ decomposition of $a_k u_k$ in $S$.  Finally, set
%Then define the \emph{rank-two exceptional sliding} $\sigma_{- \alpha,\xi_k}(p_k)$ to be
$$ q_k' = a_k' n_k' \qquad \mbox{where} \ n_k' = l_k^{-1} n_k'' l_k$$
Note that $a_k' \in A$ and $n_k' \in N_\beta^+$, so $q_k' \in H_\beta$.  Clearly $(q_k')$ is a vertical perturbation of $(q_k)$, so it is an admissible perturbation of $(p_k)$.  The conjugation by $l_k$ on $N^+_\beta$ preserves the subspace $E_1 = \lieg_{\alpha} \oplus \lieg_{\alpha + \beta}$, so $ER(q_k')$ contains $\alpha$ or $\alpha + \beta$.  If it only contains $\alpha + \beta$, then we perform a Weyl reflection $r_\beta$ to finally obtain an admissible perturbation $(q_k'')$ of $(p_k)$ with $\alpha \in ER(q_k'')$.

\bibliographystyle{amsplain}
\bibliography{karinsrefs}

\providecommand{\bysame}{\leavevmode\hbox to3em{\hrulefill}\thinspace}
\providecommand{\MR}{\relax\ifhmode\unskip\space\fi MR }
\providecommand{\MRhref}[2]{%
  \href{http://www.ams.org/mathscinet-getitem?mr=#1}{#2}
}
\providecommand{\href}[2]{#2}
\begin{thebibliography}{10}

\bibitem{lf.conf.regularity}
J.~Ferrand, \emph{Geometrical interpretation of scalar curvature and regularity
  of conformal homeomorphisms}, Differential geometry and relativity, vol.~3,
  Reidel, Dordrecht, 1976, pp.~1--44.

\bibitem{frances.lfrank1}
C.~Frances, \emph{Sur le groupe d'automorphismes des g\'eom\'etries
  paraboliques de rang 1}, Ann. Sci. \'Ecole Norm. Sup. (4) \textbf{40} (2007),
  no.~5, 741--764.

\bibitem{frances.degenerescence}
\bysame, \emph{D\'eg\'enerescence locale des transformations conformes
  pseudo-riemanniennes}, Ann. Inst. Fourier \textbf{62} (2012), no.~5,
  1627--1669.

\bibitem{fm.nilpconf}
C.~Frances and K.~Melnick, \emph{Nilpotent groups of conformal flows on compact
  pseudo-{R}iemannian manifolds}, Duke Math. J. \textbf{153} (2010), no.~3,
  511--550.

\bibitem{gromov.rgs}
M.~Gromov, \emph{Rigid transformations groups}, G\'eom\'etrie Diff\'erentielle
  (Paris, 1986) (D.~Bernard and Y.~Choquet-Bruhat, eds.), Hermann, Paris, 1988,
  pp.~65--139.

\bibitem{hochschild}
G.P. Hochschild, \emph{Basic theory of algebraic groups and {L}ie algebras},
  Graduate Texts in Mathematics, vol.~75, Springer, Berlin, 1981.

\bibitem{knapp.lie.groups}
A.W. Knapp, \emph{Lie groups: Beyond an introduction}, Birkh\"auser, Boston,
  MA, 1996.

\bibitem{kobayashi.transf}
S.~Kobayashi, \emph{Transformation groups in differential geometry}, Springer,
  Berlin, 1995.

\bibitem{lf.noncompact}
J.~Lelong-Ferrand, \emph{The action of conformal transformations on a
  {R}iemannian manifold}, Math. Ann. \textbf{304} (1996), 277--291.

\bibitem{me.frobenius}
K.~Melnick, \emph{A {F}robenius theorem for {C}artan geometries, with
  applications}, Enseign. Math. S\'er. II \textbf{57} (2011), no.~1-2, 57--89.

\bibitem{myers.steenrod}
S.B. Myers and N.~Steenrod, \emph{The group of isometries of a {R}iemannian
  manifold}, Ann. Math. \textbf{40} (1939), no.~2, 400--416.

\bibitem{nomizu.affine.transf}
K.~Nomizu, \emph{On the group of affine transformations of an affinely
  connected manifold}, Proc. Amer. Math. Soc. \textbf{4} (1953), 816--823.

\bibitem{palais.lie}
R.S. Palais, \emph{A global formulation of the {L}ie theory of transformation
  groups}, Mem. Amer. Math. Soc., vol.~22, Amer. Math. Soc., Providence, RI,
  1957.

\bibitem{ruh.finite.type}
E.A. Ruh, \emph{On the automorphism groups of a $g$-structure}, Comment. Math.
  Helv. \textbf{39} (1964), 189--264.

\bibitem{schoen.cr}
R.~Schoen, \emph{On the conformal and {CR} automorphism groups}, Geom. Funct.
  Anal. \textbf{5} (1995), no.~2, 464--481.

\bibitem{sharpe}
R.W. Sharpe, \emph{Differential geometry : {C}artan's generalization of
  {K}lein's {E}rlangen program}, Springer, New York, 1996.

\bibitem{sternberg.diffgeom}
S.~Sternberg, \emph{Lectures on differential geometry}, Prentice-Hall,
  Englewood Cliffs, NJ, 1964.

\bibitem{tanaka.parabolic}
N.~Tanaka, \emph{On the equivalence problem associated with simple graded {L}ie
  algebras}, Hokkaido Math. J. \textbf{8} (1979), 23--84.

\bibitem{cap.slovak.book.vol1}
A.~\v{C}ap and J.~Slov\'ak, \emph{Parabolic geometries {I}}, Mathematical
  Surveys and Monographs, vol. 154, American Mathematical Society, Providence,
  RI, 2009.

\end{thebibliography}

\begin{tabular}{lll}
Karin Melnick  & \quad\qquad & Charles Frances \\
Department of Mathematics & \quad\qquad & Institut de Recherche Math\'ematique Avanc\'ee  \\
4176 Campus Drive & \qquad \qquad & 7 rue Ren\'e-Descartes \\
University of Maryland & \quad\qquad & Universit\'e de Strasbourg \\
College Park, MD 20742 &\quad \qquad & 67085 Strasbourg Cedex  \\
USA &\quad \qquad & France \\
karin@math.umd.edu &\quad \qquad & frances@math.unistra.fr
\end{tabular}

\end{document}